\documentclass[11pt]{article}

\usepackage{amsthm,amsmath,amssymb,xspace,epsfig,stackrel,}
\RequirePackage[colorlinks,citecolor=blue,urlcolor=blue]{hyperref}
\theoremstyle{plain}
\newtheorem{thm}{Theorem}

\newtheorem{lam}{Lemma}
\newtheorem{col}{Corollary}

\RequirePackage[colorlinks,citecolor=blue,urlcolor=blue]{hyperref}

\usepackage[latin1]{inputenc}
\usepackage[top=3cm, bottom=3cm, left=2.5cm, right=2.5cm]{geometry}
\usepackage{amsmath, amsthm, latexsym, xspace, epsfig, stackrel}
\usepackage{amsfonts}
\usepackage{amssymb}
\usepackage{booktabs}
\usepackage{dcolumn}

\begin{document}

\title{Toward optimal model averaging in regression models with time series errors}
\author{Tzu-Chang F. Cheng, Ching-Kang Ing, \\ Shu-Hui Yu \\
\\
University of Illinois at Urbana-Champaign, \\Academia Sinica and National Taiwan University, \\ National University of Kaohsiung}
\date{}
\maketitle

\begin{abstract}
Consider a regression model with infinitely many parameters and
time series errors. We are interested in choosing weights for
averaging across generalized least squares (GLS) estimators obtained from a set of approximating models.
However, GLS estimators, depending on
the unknown inverse covariance matrix of the errors,
are usually infeasible. We therefore construct feasible generalized least squares (FGLS) estimators
using a consistent estimator of the unknown inverse matrix.
Based on this inverse covariance matrix estimator and FGLS estimators, we develop
a feasible autocovariance-corrected Mallows model averaging
criterion to select weights, thereby
providing an FGLS model averaging estimator of the true regression function.
We show that
the generalized squared error loss
of our averaging estimator is asymptotically equivalent to
the minimum one among those of GLS model averaging estimators
with the weight vectors belonging to a continuous set, which
includes the discrete weight set used in Hansen (2007)
as its proper subset.
\\
\\
\noindent {\it  JEL classification}: C22; C52
\\
\\
\noindent \textbf{KEY WORDS:}
Asymptotic efficiency;
Autocovariance-corrected Mallows model averaging;
Banded Cholesky factorization;
Feasible generalized least squares estimator;
High-dimensional covariance matrix;
Time series errors.

\end{abstract}

\def \theequation{1.\arabic{equation}}
\setcounter{equation}{0}

\section{\large Introduction}
This article is concerned with the implementation of model averaging methods in
regression models with time series errors.
We are interested in choosing weights for
averaging across generalized least squares (GLS) estimators
obtained from a set of approximating models for the true regression function.
However, GLS estimators, depending on
the unknown covariance matrix $\Sigma^{-1}_{n}$ of the errors,
are usually infeasible, where $n$ is the sample size.
We therefore construct feasible generalized least squares (FGLS) estimators
using a consistent estimator of $\Sigma^{-1}_{n}$.
Based on this inverse covariance matrix estimator and FGLS estimators, we develop
a feasible autocovariance-corrected Mallows model averaging
(FAMMA) criterion to select weights, thereby
providing an FGLS model averaging estimator of the regression function.
We show that
the generalized squared error loss
of our averaging estimator is asymptotically equivalent to
the minimum one among those of GLS model averaging estimators
with the weight vectors belonging to a continuous set, which
includes the discrete weight set used in Hansen (2007)
as its proper subset.


Let $M$ be the number of approximating models.
If the weight set only contains
standard unit vectors in $R^{M}$, then
selection of weights for model averaging
is equivalent to selection of models.
Therefore, model selection can be viewed as a special case of model averaging.
It is shown in Hansen (2007, p.1179) that when the weight set is rich enough, the optimal
model averaging estimator usually outperforms the one obtained from the optimal single model,
providing ample reason to conduct model averaging.
Another vivid example demonstrating the advantage of model averaging over model selection
is given by Yang (2007, Section 6.2.1, Figure 5).
In the case of independent errors,
asymptotic efficiency results for model selection have been reported extensively,
even when the errors are heteroskedastic or regression functions are serially correlated.
For the regression model with i.i.d. Gaussian errors, Shibata (1981) showed that
Mallows' $C_{p}$ (Mallows (1973)) and Akaike information criterion (AIC; Akaike (1974))
lead to asymptotically efficient estimators of the regression function.
By making use of Whittle's (1960) moment bounds for quadratic forms in independent variables,
Li (1987) established the asymptotic efficiency of Mallows' $C_{p}$ under much weaker assumptions on homogeneous errors.
Li's (1987) result was subsequently extended by Andrews (1991) to heteroscedastic errors.
There are also asymptotic efficiency results
established in situation where regression functions are serially correlated.
Assuming that the data are generated from an infinite order autoregressive (AR($\infty$)) process
driven by i.i.d. Gaussian noise, Shibata (1980) showed that
AIC is asymptotically efficient for independent-realization prediction.
This result was extended to non-Gaussian AR($\infty$) processes
by Lee and Karagrigoriou (2001).
Ing and Wei (2005) showed that
AIC is also asymptotically efficient for same-realization prediction.
Ing (2007) further pointed out that the same property holds for
a modification of Rissanen's accumulated
prediction error (APE, Rissanen (1986)) criterion.

Asymptotic efficiency results for model averaging
have also attracted much recent attention
from econometricians and statisticians.
Hansen (2007) proposed the Mallows model averaging (MMA)
criterion, which selects weights
for averaging across LS estimators.
Under regression models with i.i.d. explanatory vectors
and errors, he proved that the averaging estimator obtained from the MMA criterion
asymptotically attains the minimum squared error loss
among those of the LS model averaging estimators
with the weight vectors contained in a discrete set ${\cal H}_{n}(N)$ (see (2.8)),
in which $N$ is a positive integer and related to the moment restrictions of the errors.
Using the same weight set, Hansen and Racine (2012) and Liu and Okui (2013), respectively, showed that
the Jackknife model averaging (JMA) criterion and feasible HR$C_{p}$ criterion yield
asymptotically efficient LS model averaging estimators
in regression models with independent explanatory vectors and heteroscedastic errors.
Since
${\cal H}_{n}(N)$ is quite restrictive when $N$ is small,
Wan, Zhang and Zou (2010) justified MMA's asymptotic efficiency
over the continuous weight set
\begin{eqnarray}
{\cal G}_{n}=\{\boldsymbol{w}=(w_{1}, \ldots, w_{M})^{'}: w_{m} \in [0,1], \sum_{m=1}^{M}w_{m}=1\},
\end{eqnarray}
which is much more flexible than ${\cal H}_{n}(N)$.
Recently,
Ando and Li (2014) showed that
Hansen and Racine's (2012) result carries over to
high-dimensional regression models and to
a weight set more general than ${\cal G}_{n}$.

There are different types of theoretical examinations on model
averaging. Besides the approach of targeting asymptotic efficiency,
another very successful approach is minimax optimal model
combination via oracle inequalities; see, for example,
Yang (2001), Yuan and Yang (2005),
Leung and Barron (2006), and
Wang et al. (2014).

However,
all aforementioned papers, requiring the error terms to be independent,
preclude the regression model with time series errors,
which is one of the most useful models for analyzing dependent data.
In this article, we take the first step to close this gap by
introducing the FAMMA criterion
and proving its asymptotic efficiency in the sense mentioned
in the first paragraph.
However, minimax optimality results are not pursued here.
Our criterion has some distinctive features.
First, it involves estimation
of the high-dimensional inverse covariance matrix
of a stationary time series that is not directly observable.
Note that the covariance matrix of a stationary time series of length $n$
can be viewed as a high-dimensional covariance matrix because
its dimension is equivalent to the sample size.
In situations where the error process is observable
(or equivalently, the regression functions are known to be zero),
Wu and Pourahmadi (2009) proposed a banded covariance matrix estimator of
$\Sigma_{n}$ and proved its consistency under spectral norm,
which also leads to the consistency of the corresponding
inverse matrix in estimating $\Sigma^{-1}_{n}$. These results were then extended
by McMurry and Politis (2010) to tapered covariance matrix estimators.
However, since the error process is in general unobservable,
one can only estimate $\Sigma^{-1}_{n}$ (or $\Sigma_{n}$)
through the output variables.
As far as estimating $\Sigma^{-1}_{n}$ is concerned, these output variables are contaminated by unknown regression functions.
In Section 3, we propose estimating $\Sigma^{-1}_{n}$
by its banded Cholesky decomposition with the corresponding parameters
estimated {\it nonparametrically} from the {\it least squares residuals} of an increasing dimensional
approximating model. We also obtain the rate of convergence of the proposed estimator,
which plays a crucial role in proving the asymptotic efficiency of the FAMMA criterion.
Second, our criterion is justified under a continuous weight set ${\cal H}_{N}$ (see (2.10)).
While ${\cal H}_{N}$ is not as general as ${\cal G}_{n}$, as argued in Section 2, it can substantially
reduce the limitations encountered by
${\cal H}_{n}(N)$ when $N$ is small.

It is worth mentioning that
to justify MMA's asymptotic efficiency over the weight set ${\cal G}_{n}$,
Wan, Zhang and Zou (2010) required
a stringent condition on $M$; see (2.20) of Section 2.
As argued in Remark 4, this condition may preclude the approximating models
whose estimators have the minimum risk (ignoring constants).
When these models/estimators are precluded, the MMA criterion
can only select weights for a set of suboptimal models/estimators, which is obviously not desirable.
In fact, the same dilemma also arises in Ando and Li (2014), who used a similar assumption to
prove their asymptotic efficiency results.
Zhang, Wan and Zou (2013) considered model averaging problems
in regression models with dependent errors.
They adopted the JMA criterion to choose weights for a class of estimators
and showed that the criterion is asymptotically efficient
over the weight set ${\cal G}_{n}$.
Their result, however, is still reliant on a condition similar to (2.20).
In addition, the class of estimators considered in their paper,
excluding all FGLS estimators, may suffer from lack of efficiency.

The remaining paper is organized as follows.
In Section 2, we first concentrate on the case where $\Sigma_{n}$ is known.
We show in Theorem 1 that the autocovariance-corrected Mallows model averaging
(AMMA) criterion, which is the FAMMA criterion with the estimator of $\Sigma^{-1}_{n}$
replaced by $\Sigma^{-1}_{n}$ itself, is asymptotically efficient.
Since the assumptions used in Theorem 1 are rather mild,
both Corollary 2.1 of Li (1987) and Theorem 1 of Hansen (2007)
become its special case.
We then turn attention to the more practical situation
where $\Sigma_{n}$ is unknown and propose
choosing model weights by the FAMMA criterion.
It is shown in Theorem 2 of Section 2 that
the FAMMA criterion is asymptotically efficient as long as the corresponding estimator of $\Sigma^{-1}_{n}$
has a sufficiently fast convergence rate.
In Section 3, we provide a consistent estimator of $\Sigma^{-1}_{n}$
based on its banded Cholesky decomposition, and derive the estimator's convergence rate under various situations.
In Section 4, the asymptotic efficiency of the FAMMA criterion with $\Sigma^{-1}_{n}$
estimated by the method proposed in Section 3 is established.
Finally, we conclude in Section 5.
All proofs are relegated to the Appendix in order to maintain the flow of exposition.


\def \theequation{2.\arabic{equation}}
\setcounter{equation}{0}

\section{\large The AMMA and FAMMA criteria}

\indent
Consider a regression model with infinitely many parameters,
\begin{align}
y_{t}=\sum^{\infty}_{j=1}\theta_{j}x_{tj}+e_{t}=\mu_{t}+e_{t}, t=1, \ldots, n,
\end{align}
where
$\mu_{t}=\sum^{\infty}_{j=1}\theta_{j}x_{tj}$,
$\boldsymbol{x_{t}}=(x_{t1},x_{t2}, \ldots)^{'}$ is the explanatory vector
with $\sup_{t\geq 1, j\geq 1}\mathrm{E}(x_{tj}^{2})<\infty$,
$\theta_{j}, j\geq 1$ are unknown parameters
satisfying $\sum_{j=1}^{\infty}|\theta_{j}|<\infty$,
and
$\{e_{t}\}$,
independent of $\{\boldsymbol{x_{t}}\}$,
is an unobservable stationary process with zero mean and finite variance.
In matrix notation, $Y_{n}=\boldsymbol{\mu}_{n}+\textbf{e}_{n}$,
where $Y_{n}=(y_{1},...,y_{n})'$, $\boldsymbol{\mu}_{n}=(\mu_{1},...,\mu_{n})'$, and $\textbf{e}_{n}=(e_{1},...,e_{n})'$.
The central focus of this paper is to explore how and to what extent the model averaging
can be implemented in the presence of time series errors.

Let $m=1, \ldots, M$
be a set of approximating models of (2.1),
where the $m$th model uses the first $k_{m}$ elements of $\{\boldsymbol{x}_{t}\}$
with $1\leq k_{1}<k_{2}<\cdots <k_{M}<n$ and
$M$ is allowed to grow to infinity with the sample size $n$.
Assume that $\Sigma_{n}=\mathrm{E}(\textbf{e}_{n}\textbf{e}'_{n})$ is known and $\Sigma^{-1}_{n}$ exists.
Then the generalized least squares (GLS) estimator of the regression coefficient vector
in the $m$th approximating model is given by
$\hat{\Theta}^{*}_{m}=(X_{m}'\Sigma_{n}^{-1}X_{m})^{-1}X_{m}'\Sigma^{-1}_{n}Y_{n}$,
and the resultant estimate of $\boldsymbol{\mu}_{n}$ is
$\hat{\boldsymbol{\mu}}_{n}(m)=P^{*}_{m}Y_{n}$,
where
$X_{m}=(x_{ij})_{1\leq i\leq n, 1\leq j \leq k_{m}}$,
$P^{*}_{m}=X_{m}(X_{m}'\Sigma_{n}^{-1}X_{m})^{-1}X_{m}'\Sigma^{-1}_{n}$, and
$X_{M}$ is assumed to be almost surely (a.s.) full rank throughout the paper.
The model averaging estimator of $\boldsymbol{\mu}_{n}$
based on the $M$th approximating models is $\hat{\boldsymbol{\mu}}_{n}(\boldsymbol{w})=P^{*}(\boldsymbol{w})Y_{n}$,
where
$\boldsymbol{w} \in {\cal G}_{n}$
and $P^{*}(\boldsymbol{w})=\sum_{m=1}^{M}w_{m}P^{*}_{m}$.
To evaluate the performance of $\hat{\boldsymbol{\mu}}_{n}(\boldsymbol{w})$,
we use the generalized squared error (GSE) loss
\begin{eqnarray*}
L^{*}_{n}(\boldsymbol{w})=(\hat{\boldsymbol{\mu}}_{n}(\boldsymbol{w})-\boldsymbol{\mu}_{n})^{'}\Sigma^{-1}_{n}
(\hat{\boldsymbol{\mu}}_{n}(\boldsymbol{w})-\boldsymbol{\mu}_{n}).
\end{eqnarray*}
This loss function is a natural generalization of
Hansen's (2007) average squared error
in the sense that
$L^{*}_{n}(\boldsymbol{w})$
reduces to the latter when
$\Sigma^{-1}_{n}$
is replaced by the $n \times n$
identity matrix.
Through this generalization,
it is easy to establish a connection
between our results and
some classical asymptotic efficiency results
on model averaging/selection,
thereby leading to a more comprehensive understanding of this research field.
For further discussion, see Remarks 1 and 2 below.
On the other hand, when the future values of $y_{t}$
are entertained instead of the regression function $\boldsymbol{\mu}_{n}$,
Wei and Yang (2012) proposed several different loss functions from a prediction point of view.
Compared with the squared errors,
their loss functions are particularly suitable for dealing with outliers.

The next lemma provides a representation for the conditional risk,
\begin{eqnarray*}
R^{*}_{n}(\boldsymbol{w})=\mathrm{E}(L^{*}_{n}(\boldsymbol{w})|\boldsymbol{x}_{1}, \ldots, \boldsymbol{x}_{n})\equiv
\mathrm{E}_{\boldsymbol{x}}(L^{*}_{n}(\boldsymbol{w})),
\end{eqnarray*}
which is an extension of Lemma 2 of Hansen (2007)
to the case of dependent errors.
\begin{lam}
Assume {\rm (2.1)} and $\Sigma^{-1}_{n}$ exists. Then, for any
$\boldsymbol{w} \in {\cal G}_{n}$,
\begin{eqnarray}
R^{*}_{n}(\boldsymbol{w})=\sum_{m=1}^{M}\sum_{l=1}^{M}w_{m}w_{l}
[\boldsymbol{\mu}^{'}_{n}\Sigma^{-1/2}_{n}(I-P_{\max\{m, l\}})\Sigma^{-1/2}_{n}\boldsymbol{\mu}_{n}+\min\{k_{m}, k_{l}\}],
\end{eqnarray}
where $P_{j}=\Sigma^{-1/2}_{n}X_{j}(X^{'}_{j}\Sigma^{-1}_{n}X_{j})^{-1}X^{'}_{j}\Sigma^{-1/2}_{n}$
is the orthogonal projection matrix for the column space of $\Sigma^{-1/2}_{n}X_{j}$.
\end{lam}

To choose a data-driven weight vector
asymptotically minimizing $L^{*}_{n}(\boldsymbol{w})$ over a suitable weight set ${\cal H} \subseteq {\cal G}_{n}$,
we propose the AMMA criterion,
\begin{eqnarray}
C^{*}_{n}(\boldsymbol{w})=(Y_{n}-\hat{\boldsymbol{\mu}}_{n}(\boldsymbol{w}))^{'}\Sigma^{-1}_{n}(Y_{n}-\hat{\boldsymbol{\mu}}_{n}(\boldsymbol{w}))+
2\sum_{m=1}^{M}w_{m}k_{m}.
\end{eqnarray}
Note that the AMMA criterion is a special case of the criterion given in (6.2$^{*}$) of Andrews (1991) with $M_{n}(h)=P^{*}(\boldsymbol{w})$
and $W=\Sigma^{-1}_{n}$.
It reduces to the MMA criterion when
\begin{eqnarray}
\{e_{t}\} \,\,\mbox{is a sequence of i.i.d. random variables with} \,\,\mathrm{E}(e_{1})=0, 0<\mathrm{E}(e^{2}_{1})=\sigma^{2}<\infty.
\end{eqnarray}
Recently, Liu, Okui and Yoshimura (2013) also suggested using $C^{*}_{n}(\boldsymbol{w})$
to choose weight vectors in situations where $e_{t}$ are independent but possibly heteroscedastic.
While the AMMA criterion is not new to the literature,
the question of whether its minimizer can (asymptotically) minimize $L^{*}_{n}(\boldsymbol{w})$
seems rarely discussed, in particular when ${\cal H}$ is uncountable.


\vspace{0.2cm}
Recall that by assuming (2.4),
\begin{eqnarray}
\xi_{n}=\inf_{\boldsymbol{w} \in {\cal G}_{n}} R^{*}_{n}(\boldsymbol{w}) \to \infty \,\,\mbox{a.s.},
\end{eqnarray}
and
\begin{eqnarray}
\mathrm{E}(|e_{1}|^{4(N+1)}|\boldsymbol{x}_{1})< \kappa<\infty \,\,\mbox{a.s.}, \,\,\mbox{for some positive integer}\,\, N,
\end{eqnarray}
Hansen (2007, Theorem 1) showed that
$C^{*}_{n}(\boldsymbol{w})$
is asymptotically efficient in the sense that
\begin{eqnarray}
\frac{L^{*}_{n}(\bar{\boldsymbol{w}}_{n})}{\inf_{\boldsymbol{w} \in {\cal H}_{n}(N)}L^{*}_{n}(\boldsymbol{w})} \to_{p} 1,
\end{eqnarray}
where $\to_{p}$ denotes convergence in probability,
\begin{eqnarray}
{\cal H}_{n}(N)=\{\boldsymbol{w}: w_{m} \in \{0, 1/N, 2/N, \ldots 1\}, \sum_{m=1}^{M}w_{m}=1\},
\end{eqnarray}
and
\begin{eqnarray*}
\bar{\boldsymbol{w}}_{n}=\arg \min_{\boldsymbol{w} \in {\cal H}_{n}(N)} C^{*}_{n}(\boldsymbol{w}).
\end{eqnarray*}
Equation (2.7) gives a positive answer to the above question
in the special case where $\Sigma_{n}=\sigma^{2}I_{n}$
and ${\cal H}={\cal H}_{n}(N)$ is a discrete set.
When (2.6) holds for sufficiently large $N$, the restriction of ${\cal G}_{n}$
to ${\cal H}_{n}(N)$
is not an issue of overriding concern because the grid points $i/N, i=0, \ldots, N$,
in ${\cal H}_{n}(N)$ is dense
enough to provide a good approximation for the optimal
weight vector among
\begin{eqnarray}
\bar{{\cal H}}_{n}(N)=\{\boldsymbol{w}: w_{m} \in [0,1], 1\leq \sum_{m=1}^{M}I_{\{w_{m}\neq 0\}} \leq N, \sum_{m=1}^{M}w_{m}=1\},
\end{eqnarray}
and hence among ${\cal G}_{n}$.
We call $\bar{{\cal H}}_{n}(N)$ {\it continuous extension} of ${\cal H}_{n}(N)$
because it satisfies ${\cal H}_{n}(N) \subseteq \bar{{\cal H}}_{n}(N)$
and $a\boldsymbol{w}_{1}+b\boldsymbol{w}_{2} \in \bar{{\cal H}}_{n}(N)$,
for any
$0\leq a, b \leq 1$ with $a+b=1$
and any
$\boldsymbol{w}_{1}=(w_{11}, \ldots w_{M1})^{'}$ and $\boldsymbol{w}_{2}=
(w_{12}, \ldots w_{M2})^{'} \in \bar{{\cal H}}_{n}(N)$ with $\sum_{m=1}^{M}
|I_{\{w_{m1}\neq 0\}}-I_{\{w_{m2}\neq 0\}}|=0$.
It is shown in Hansen (2007, p.1179) that even when $N=2$,
the optimal weight vector in $\bar{{\cal H}}_{n}(N)$
can yield an averaging estimator outperforming the one based on the optimal single model, except in some special cases.

On the other hand, when (2.6) holds only for moderate or small
values of $N$, not only the optimal weight vector in ${\cal G}_{n}$
but also that in $\bar{{\cal H}}_{n}(N)$ cannot be well approximated
by the elements in ${\cal H}_{n}(N)$.
As a result, the advantage of model averaging over model selection
becomes less apparent. To rectify this deficiency, we introduce the
following continuous extension of ${\cal H}_{n}(N)$,
\begin{eqnarray}
{\cal H}_{N}=\bigcup_{l=1}^{N}{\cal H}_{(l)},
\end{eqnarray}
where
\begin{eqnarray*}
{\cal H}_{(l)}=\{\boldsymbol{w}:
\underline{\delta}\leq w_{i}I_{\{w_{i}\neq 0\}}\leq 1, \sum_{i=1}^{M}I_{\{w_{i} \neq 0\}}=l, \sum_{m=1}^{M}w_{m}=1\},
\end{eqnarray*}
with $0<\underline{\delta} < 1/N$.
The number of non-zero component is $1\leq l \leq N$
for any vector in ${\cal H}_{(l)}$.
Therefore, this weight set leads to sparse combinations of $\hat{\boldsymbol{\mu}}_{n}(m), m=1, \ldots M$.
For a detailed discussion on
sparse combinations from the minimax viewpoint,
see Wang et al. (2014).
We will show in Theorem 1 that
\begin{eqnarray}
\frac{L^{*}_{n}(\tilde{\boldsymbol{w}}_{n})}{\inf_{\boldsymbol{w} \in {\cal H}_{N}}L^{*}_{n}(\boldsymbol{w})} \to_{p} 1,
\end{eqnarray}
without the restriction $\Sigma_{n}=\sigma^{2}I_{n}$, where
\begin{eqnarray*}
\tilde{\boldsymbol{w}}_{n}=\arg \inf_{\boldsymbol{w} \in {\cal H}_{N}} C^{*}_{n}(\boldsymbol{w}).
\end{eqnarray*}

It is important to be aware that ${\cal H}_{N}$ can
inherit the benefits of $\bar{{\cal H}}_{n}(N)$
mentioned previously because
the difference between the two sets
can be made arbitrarily small by making
$\underline{\delta}$ sufficiently close to 0.
Technically speaking,
a nonzero (regardless of how small) $\underline{\delta}$ enables us to
establish some sharp uniform probability bounds
through replacing $R^{*}_{n}(\boldsymbol{w})$
by suitable model selection risks (see (A.7)),
thereby overcoming the difficulties arising from the uncountablility of
${\cal H}_{N}$.
In Remarks 3 and 4 after Theorem 1,
we will also discuss the asymptotic efficiency of $C^{*}_{n}(\boldsymbol{w})$ over more general weight sets
such as ${\cal G}_{n}$ and its variants.
The following assumptions on $\{e_{t}\}$ are needed in our analysis.
As shown in the Appendix,
these assumptions allow us to derive
sharp bounds for the moments of quadratic forms in $\{e_{t}\}$
using the first moment bound theorem of Findley and Wei (1993).

\vspace{0.3cm}
\noindent
{\bf Assumption 1.}
$\{e_{t}\}$ is a sequence of stationary time series
with autocovariance function (ACF) $\gamma_{j}=\mathrm{E}(e_{t}e_{t+j})$
satisfying $\sum_{j=-\infty}^{\infty}\gamma^{2}_{j}<\infty$, and admits a linear representation
\begin{align}
e_{t}=\alpha_{t}+\sum^{\infty}_{k=1}\beta_{k}\alpha_{t-k}
\end{align}
in terms of the ${\cal F}_{t}$-measurable random variables $\alpha_{t}$,
where ${\cal F}_{t}, -\infty<t<\infty$ is an increasing sequence
of $\sigma$-fields of events. Moreover,
$\{\alpha_{t}\}$ satisfies the following properties with probability 1:
\begin{description}
\item(M1) $E(\alpha_{t}|{\cal F}_{t-1})=0$.
\item(M2) $E(\alpha^{2}_{t}|{\cal F}_{t-1})=\sigma^{2}_{\alpha}$.
\item(M3) There exist a positive integer $N$ and a positive number $S>4N$ such that
for some constant $0<C_{S}<\infty$,
\begin{align}
&\sup_{-\infty < t<\infty}E(|\alpha_{t}|^{S}|{\cal F}_{t-1})\leq C_{S}.
\end{align}
\end{description}

\noindent
{\bf Assumption 2.}
The spectral density function of $\{e_{t}\}$,
\begin{eqnarray}
f_{e} (\lambda)=(\sigma^{2}_{\alpha}/2\pi)
|\sum_{j=0}^{\infty}\beta_{j}e^{-ij\lambda}|^{2} \neq 0
\end{eqnarray}
for all $-\pi < \lambda\leq \pi$, where $\beta_{0}=1$.
Moreover,
\begin{eqnarray}
\sum_{j=0}^{\infty}|\beta_{j}|<\infty.
\end{eqnarray}

\vspace{0.3cm}
We are now in a position to state Theorem 1.
\begin{thm}
Assume {\rm Assumptions 1 and 2} in which $N$ in (M3) is fixed.
Let
\begin{eqnarray*}
D_{n}(m)=\boldsymbol{\mu}^{'}_{n}\Sigma^{-1/2}_{n}(I-P_{m})\Sigma^{-1/2}_{n}\boldsymbol{\mu}_{n}+k_{m} \,\,
{\rm and}\,\, k^{*}_{n}=\min_{1\leq m \leq M} D_{n}(m).
\end{eqnarray*}
Suppose
\begin{eqnarray}
 k^{*}_{n} \to \infty \,\,{\rm a.s.}
\end{eqnarray}
Then, {\rm (2.11)} follows.
\end{thm}

\vspace{0.2cm}
A few comments on Theorem 1 are in order.
\\
\\
{\bf Remark 1.}
Theorem 1 generalizes Theorem 1 of Hansen (2007) in several directions.
First, (2.4) is a special case of (2.12), with $\beta_{k}=0$ for all $k \geq 1$.
Second, the discrete weight set ${\cal H}_{n}(N)$ is extended to its continuous extension ${\cal H}_{N}$.
Third, when (2.4) holds and $\{e_{t}\}$ is independent of $\{\boldsymbol{x}_{t}\}$,
the moment condition (2.13) is milder than (2.6).
Fourth , (2.5) is weakened to (2.16), which is much easier to verify.
Note that $k^{*}_{n}$ can be viewed as
an index of the amount of information contained in the candidate models.
Therefore, (2.16) is quite natural from
the estimation theoretical viewpoint; see, e.g., Lai and Wei (1982), Yu, Lin and Cheng (2012) and Chan, Huang and Ing (2013).
Suppose there exists a non-random and non-negative function $Q(m)$ satisfying
\begin{eqnarray*}
\sup_{1\leq m \leq M}\left|\frac{\boldsymbol{\mu}^{'}_{n}\Sigma^{-1/2}_{n}(I-P_{m})\Sigma^{-1/2}_{n}\boldsymbol{\mu}_{n}}{n}-Q(m)\right|
\to 0, \,\,\mbox{a.s.}
\end{eqnarray*}
Then, (2.16) is fulfilled if
$Q(m) \neq 0$ for all $m$ and $\lim_{m \to \infty}Q(m)=0$,
which essentially require that all candidate models are misspecified, but
those which have many parameters
can give good approximations of the true model.
\\
\\
{\bf Remark 2.} Corollary 2.1 of Li (1987) also becomes a special case of Theorem 1.
To see this, note that under (2.4), (2.16) and
\begin{eqnarray}
\mathrm{E}(e^{8}_{1})<\infty,
\end{eqnarray}
Li's (1987) Corollary 2.1 shows that (2.11) holds with $N=1$, namely,
$C^{*}_{n}(\boldsymbol{w})$ is asymptotically efficient for model selection.
However, since (2.4) and (2.17) also imply Assumptions 1 and 2,
Li's conclusion readily follows from Theorem 1.
\\
\\
{\bf Remark 3.} It is far from being trivial to extend (2.11) to
\begin{eqnarray}
\frac{L^{*}_{n}(\boldsymbol{w}^{0})}{\inf_{\boldsymbol{w} \in {\cal G}_{n}}L^{*}_{n}(\boldsymbol{w})} \to_{p} 1,
\end{eqnarray}
where
\begin{eqnarray*}
\boldsymbol{w}^{0}=\arg \inf_{\boldsymbol{w} \in {\cal G}_{n}} C^{*}_{n}(\boldsymbol{w}).
\end{eqnarray*}
Alternatively, if $\alpha_{t}$ have light-tailed distributions, such as those
described in (C2) and (C3) of Ing and Lai (2011), then by using the exponential probability inequalities developed
in the same papers in place of the moment inequalities given in the proof of Theorem 1,
it can be shown that (2.11) holds with $\underline{\delta}$ tending to 0 and
$N=M$ tending to $\infty$ sufficiently slowly with $n$.
The details, however, are not reported here due to space constraints.
\\
\\
{\bf Remark 4.}
When (2.4) holds true,
(2.18) has been developed in Theorem 1' of Wan, Zhang and Zou (2010) under
\begin{eqnarray}
\mathrm{E}(|e|^{4G}_{t}|\boldsymbol{x}_{t})\leq \kappa <\infty \,\,\mbox{a.s.},
\end{eqnarray}
for some integer $1\leq G<\infty$, and
\begin{eqnarray}
M\xi^{-2G}_{n}\sum_{m=1}^{M}D_{n}^{G}(m) \to 0 \,\,\mbox{a.s.}
\end{eqnarray}
Unfortunately, (2.20), imposing a stringent restriction on $M$,
often precludes models having small GSE losses.
To see this, assume that $\boldsymbol{x}_{t}$ are nonrandom and
the $m$th approximating model contains the first $m$ regressors, namely $k_{m}=m$.
Assume also that
\begin{eqnarray}
D_{n}(m)=nm^{-a}+m
\end{eqnarray}
and the $G$ in (2.19) is greater than $1/a$ for some $a \geq 1$.
It is easy to show that $D_{n}(m)$ is minimized by $m \sim (an)^{1/(1+a)}$,
yielding the optimal rate of $D_{n}(m)$, $n^{1/(1+a)}$.
Moreover,
as will be clear from (A.5), $D_{n}(m)=R^{*}_{n}(\boldsymbol{v}_{m})$
is asymptotically equivalent to $L^{*}_{n}(\boldsymbol{v}_{m})$,
where $\boldsymbol{v}_{m}$ is the $m$th standard unit vector in $R^{M}$.
Hence the optimal rate of $L^{*}_{n}(\boldsymbol{v}_{m})$
is also $n^{1/(1+a)}$, which
is achievable by any approximating model whose number of regressors
$m$ satisfying
\begin{eqnarray}
c_{1}n^{1/(1+a)}\leq m \leq c_{2}n^{1/(1+a)}, \,\,\mbox{for some}\,\,0<c_{1}<c_{2}<\infty.
\end{eqnarray}

If $M \sim c_{0}n^{1/(1+a)}$, where $c_{0}$ is any positive number,
then for $G>1/a$ with $a\geq 1$, there exists $c_{3}>0$ such that
$M\xi^{-2G}_{n}\sum_{m=1}^{M}D_{n}^{G}(m)\geq c_{3}Mn^{G}/n^{2G/(1+a)} \to \infty$ as $n \to \infty$,
which violates (2.20). In fact, it is shown in Example 2 of
Wan, Zhang and Zou (2010) that a sufficient condition for (2.20) to hold is
$M=O(n^{v})$ with $v<G/(1+2aG)<1/(1+a)$.
These facts reveal that
all models with
$L^{*}_{n}(\boldsymbol{v}_{m})$ achieving the optimal rate $n^{1/(1+a)}$ (or equivalently, with
$m$ obeying (2.22)) are excluded by (2.20).
Under such a situation,
$C^{*}_{n}(\boldsymbol{w})$ can only select weights for a set of
suboptimal models.
Therefore, it is hard to conclude from their Theorem 1' that
$\hat{\boldsymbol{\mu}}_{n}(\boldsymbol{w}^{0})$'s
GSE loss is asymptotically smaller than that of $\hat{\boldsymbol{\mu}}_{n}(m)$ with $m$ satisfying (2.22),
even though this theorem guarantees $C^{*}_{n}(\boldsymbol{w})$'s asymptotic efficiency in the sense of (2.18).
On the other hand, since Theorem 1 does not impose any restrictions similar to (2.20),
one is free to choose $M \sim \bar{C}n^{1/(1+a)}$, with $\bar{C}$ sufficiently large,
so as to include the optimal model $m \sim (an)^{1/(1+a)}$. In addition, by noticing
$k^{*}_{n}=\min_{1\leq m \leq M}D_{n}(m) \to \infty$ as $n \to \infty$,
we know from Theorem 1 that
$\tilde{\boldsymbol{w}}_{n}$ satisfies (2.11),
and hence $\hat{\boldsymbol{\mu}}_{n}(\tilde{\boldsymbol{w}}_{n})$
asymptotically outperforms the best one among
$\hat{\boldsymbol{\mu}}_{n}(m), 1\leq m \leq \bar{C}n^{1/(1+a)}$,
in terms of GSE loss.
When $a$ is unknown, the bound $M \sim \bar{C}n^{1/(1+a)}$ is infeasible.
However, if a strict lower bound for $a$, say $\underline{a}$, is known a priori,
then the same conclusion still holds for $M\sim \bar{C}^{*}n^{1/(1+\underline{a})}$
with any $\bar{C}^{*}>0$.
\\
\\
{\bf Remark 5.}
It is worth noting that
estimating the weight that minimizes
$L^{*}_{n}(\boldsymbol{w})$
will generally introduce a variance inflation factor,
which may prevent us from obtaining the asymptotic efficiency.
Under independent errors, a recent paper by Wang et al (2014) gives a
comprehensive discussion of this matter from the minimax viewpoint.
In fact, pursuing the minimax optimal rate
is more relevant than the asymptotic efficiency
in the presence of a large variance inflation factor.
On the other hand, one can still attain the asymptotic efficiency
by substantially suppressing this factor through:
(i) reducing the size of the weight set and (ii) reducing the number of the candidate variables,
which have been taken by Hansen (2007) and Wan et al. (2010), respectively.
Unfortunately, the limitations imposed
on the size of the weight set or $M$
by these authors
are too stringent, and hence may lead to suboptimal results,
as discussed previously.
Theorem 1 takes the first approach and provides
a somewhat striking result that asymptotically efficient
model averaging is still achievable
under a continuous/uncountable weight set,
which is in sharp contrast to Hansen's (2007) discrete/countable weight set.
The theoretical underpinnings of Theorem 1 are some sharp uniform probability bounds,
which are presented in the Appendix
and established based on a mild lower bound condition on the weight set described in (2.10).

\vspace{0.2cm}
In the case where $\Sigma_{n}$ is unknown, the asymptotic efficiency
of $C^{*}_{n}(\boldsymbol{w})$ developed in Theorem 1 becomes practically irrelevant.
However, if there exists a consistent estimate, $\hat{\Sigma}^{-1}_{n}$, of $\Sigma^{-1}_{n}$,
then the corresponding FGLS estimator
of $\boldsymbol{\mu}$ based on the $m$th approximating model
is $\hat{P}^{*}_{m}Y_{n}$, where
$\hat{P}^{*}_{m}=X_{m}(X^{'}_{m}\hat{\Sigma}^{-1}_{n}X_{m})^{-1}X^{'}_{m}\hat{\Sigma}^{-1}_{n}$.
Moreover, the FAMMA criterion,
\begin{eqnarray}
\hat{C}^{*}_{n}(\boldsymbol{w})=(Y_{n}-\hat{\boldsymbol{\mu}}^{*}_{n}(\boldsymbol{w}))^{'}
\hat{\Sigma}^{-1}_{n}(Y_{n}-\hat{\boldsymbol{\mu}}^{*}_{n}(\boldsymbol{w}))+
2\sum_{m=1}^{M}w_{m}k_{m},
\end{eqnarray}
can be used in place of $C^{*}_{n}(\boldsymbol{w})$
to perform model averaging, where
\begin{eqnarray*}
\hat{\boldsymbol{\mu}}^{*}_{n}(\boldsymbol{w})=\hat{P}^{*}(\boldsymbol{w})Y_{n}=\sum_{m=1}^{M}w_{m}\hat{P}^{*}_{m}Y_{n}
\end{eqnarray*}
is the FGLS model averaging estimator of $\boldsymbol{\mu}_{n}$ given $\boldsymbol{w}$.
Define
\begin{eqnarray*}
L^{F}_{n}(\boldsymbol{w})=(\hat{\boldsymbol{\mu}}^{*}_{n}(\boldsymbol{w})-\boldsymbol{\mu}_{n})^{'}\Sigma^{-1}_{n}
(\hat{\boldsymbol{\mu}}^{*}_{n}(\boldsymbol{w})-\boldsymbol{\mu}_{n}),
\end{eqnarray*}
and
\begin{eqnarray*}
\hat{\boldsymbol{w}}_{n}=\arg \inf_{\boldsymbol{w} \in {\cal H}_{N}}\hat{C}^{*}_{n}(\boldsymbol{w}).
\end{eqnarray*}

In the next theorem,
we shall show that
as long as
$\hat{\Sigma}^{-1}_{n}$ converges to
$\Sigma^{-1}_{n}$ sufficiently fast in terms of spectral norm,
(2.11) still holds with
$L^{*}_{n}(\tilde{\boldsymbol{w}}_{n})$
replaced by
$L^{F}_{n}(\hat{\boldsymbol{w}}_{n})$.

\begin{thm}
Assume that
{\rm Assumptions 1 and 2} hold and
there exists a sequence of positive numbers $\{b_{n}\}$
satisfying $b_{n}=o(n^{1/2})$
such that
\begin{eqnarray}
n\|\hat{\Sigma}^{-1}_{n}-\Sigma^{-1}_{n}\|^{2}=O_{p}(b^{2}_{n}),
\end{eqnarray}
where for the $p \times p$ matrix $A$,
$\|A\|^{2}=\sup_{\boldsymbol{z} \in R^{p}, \|\boldsymbol{z}\|=1}\boldsymbol{z}^{'}A^{'}A\boldsymbol{z}$
with $\|\boldsymbol{z}\|$ denoting the Euclidean norm of $\boldsymbol{z}$.
Moreover, suppose that there exists $2N/S<\theta<1/2$ such that
\begin{eqnarray}
\lim_{n \to \infty}\frac{b^{2}_{n}}{k^{*^{1-2\theta}}_{n}}= 0 \,\,\mbox{a.s.}
\end{eqnarray}
Then,
\begin{eqnarray}
\frac{L^{F}_{n}(\hat{\boldsymbol{w}}_{n})}{\inf_{\boldsymbol{w} \in {\cal H}_{N}}L^{*}_{n}(\boldsymbol{w})} \to_{p} 1.
\end{eqnarray}
\end{thm}

\vspace{0.2cm}
Below are some comments regarding Theorem 2.
\\
\\
{\bf Remark 6.} In the next section, (2.24) will be established for $\hat{\Sigma}^{-1}_{n}=\hat{\Sigma}^{-1}_{n}(q_{n})$,
where $\hat{\Sigma}^{-1}_{n}(q_{n})$, defined in (3.4),
is obtained by the $q_{n}$-banded Cholesky decomposition of $\Sigma^{-1}_{n}$
with the parameters in the Cholesky factors estimated nonparametrically from the least squares residuals of an increasing dimensional
approximating model. As will be seen later, the order of the magnitude of $b_{n}$
associated with $\|\hat{\Sigma}^{-1}_{n}(q_{n})-\Sigma^{-1}_{n}\|$
 can vary depending on the strength of the dependence
of $\{e_{t}\}$.
\\
\\
{\bf Remark 7.}
Zhang, Wan and Zou (2013)
considered the model averaging estimator
$\tilde{\boldsymbol{\mu}}_{n}(\boldsymbol{w})=\sum_{m=1}^{M}w_{m}\tilde{\boldsymbol{\mu}}_{n}(m)$
of $\boldsymbol{\mu}$, where
\begin{eqnarray}
\tilde{\boldsymbol{\mu}}_{n}(m)=P^{X}_{m}Y_{n}
\end{eqnarray}
the estimator corresponding to the $m$th approximating model
and $P^{X}_{m}$, an $n \times n$ matrix, depends on $\{\boldsymbol{x}_{t}\}$ only.
They evaluated
the performance of $\tilde{\boldsymbol{\mu}}_{n}(\boldsymbol{w})$ using
the usual squared error loss,
\begin{eqnarray*}
L_{n}(\boldsymbol{w})=\|\tilde{\boldsymbol{\mu}}_{n}(\boldsymbol{w})-\boldsymbol{\mu}\|^{2},
\end{eqnarray*}
and showed in their Theorem 2.1 that
\begin{eqnarray}
\frac{L_{n}(\hat{\boldsymbol{w}}^{(J)}_{n})}{\inf_{\boldsymbol{w} \in {\cal G}_{n}}L_{n}(\boldsymbol{w})} \to_{p} 1,
\end{eqnarray}
where $\hat{\boldsymbol{w}}^{(J)}_{n}$ is obtained from the JMA criterion
(defined in equation (4) of their paper).
While the weight set in (2.28) is more general than that in (2.26),
an assumption similar to (2.20) is required in their proof of (2.28).
In addition,
(2.27), excluding all FGLS estimators (since
$\hat{\Sigma}^{-1}_{n}$ depends on both $\{\boldsymbol{x}_{t}\}$
and $Y_{n}$), can suffer from lack of efficiency in estimating $\boldsymbol{\mu}$.
\\
\\
{\bf Remark 8.}
When $e_{t}$ are independent random variables with $\mathrm{E}(e_{t})=0$ for all $t$
and possibly unequal $\mathrm{E}(e^{2}_{t})=\sigma^{2}_{t}$,
Liu, Okui and Yoshimura (2013, Theorem 4) obtained a weaker version of (2.26),
\begin{eqnarray*}
\frac{L^{F}_{n}(\hat{\boldsymbol{w}}_{n})}{\inf_{\boldsymbol{w} \in {\cal H}_{n}(N)}L^{*}_{n}(\boldsymbol{w})} \to_{p} 1,
\end{eqnarray*}
in which $\hat{\Sigma}^{-1}_{n}=\mbox{diag}(\hat{\sigma}^{-2}_{1}, \cdots, \hat{\sigma}^{-2}_{n})$
with $\hat{\sigma}^{-2}_{t}$ satisfying
\begin{eqnarray}
\sup_{1\leq t\leq n}(\hat{\sigma}^{-2}_{t}-\sigma^{-2}_{t})^{2} =O_{p}(n^{-1}),
\end{eqnarray}
among other conditions.
However, since $n^{-1}$ is a parametric rate,
certain parametric assumptions on $\sigma^{2}_{t}, 1\leq t \leq n$, are required to ensure
(2.29). In addition, their proof, relying crucially on Theorem 2 of Whittle (1960),
is not directly applicable to dependent data.
\\
\\
{\bf Remark 9.}
Assumption (2.25) is a strengthened version of (2.16).
It essentially says that
the (normalized) estimation error of $\hat{\Sigma}^{-1}_{n}$ must be dominated by
the amount of information contained in the candidate models in a certain way.
This type of assumption seems indispensable for the FAMMA criterion to preserve the features of its infeasible counterpart.
\\
\\
{\bf Remark 10.}
Throughout this paper,
the only assumption that we impose on $\{\boldsymbol{x}_{t}\}$
is $\sup_{t \geq 1, j\geq 1}\mathrm{E}(|x_{tj}|^{\nu}) \\ <\infty$ for some $2 \leq \nu <\infty$,
in addition to the (a.s.) nonsingularity of $X_{M}$.
Therefore, $\{\boldsymbol{x}_{t}\}$
can be nonrandom, serially independent or serially dependent.

\def \theequation{3.\arabic{equation}}
\setcounter{equation}{0}

\section{A consistent estimate of $\Sigma^{-1}_{n}$ based on the Cholesky decomposition.}
In this section, we shall construct a consistent estimator of $\Sigma^{-1}_{n}$
based on its banded Cholesky decomposition.
Note first that according to (2.12), (2.14) and (2.15),
$e_{t}$ has an AR($\infty$) representation,
\begin{eqnarray}
\sum_{j=0}^{\infty}a_{j}e_{t-j}=\alpha_{t},
\end{eqnarray}
where $a_{0}=1$, $\sum_{j=0}^{\infty}a_{j}z^{j}=(\sum_{j=0}^{\infty}\beta_{j}z^{j})^{-1} \neq 0$ for all $|z|\leq 1$
and $\sum_{j=0}^{\infty}|a_{j}|<\infty$; see Zygmund (1959).
If an AR($k$), $k \geq 1$, model is used to approximate model (3.1), then
the corresponding best (in the sense of mean squared error) AR coefficients
are given by $-(a_{1}(k),\ldots,a_{k}(k))^{'}$, where
\begin{equation*}
(a_{1}(k),\ldots,a_{k}(k))^{'}=\text{arg$\min$}_{(c_{1},\ldots,c_{k})^{'}\in
R^{k}}\mathrm{E}(e_{t}+c_{1}e_{t-1}+\cdots+c_{k}e_{t-k})^{2}.
\end{equation*}
Define
$\sigma_{k}^{2}=\mathrm{E}(e_{t}+a_{1}(k)e_{t-1}+\cdots+a_{k}(k)e_{t-k})^{2}$.
Then, the modified Cholesky decomposition for $\Sigma^{-1}_{n}$ is
\begin{eqnarray}
\Sigma_{n}^{-1}=\mathbf{T}_{n}'\mathbf{D}_{n}^{-1}\mathbf{T}_{n},
\end{eqnarray}
where
$$\mathbf{D}_{n}=diag(\gamma_{0}, \sigma_{1}^{2}, \sigma_{2}^{2}, \cdots, \sigma_{n-1}^{2}),$$
and
$\mathbf{T}_{n}=(t_{ij})_{1\leq i,j\leq n}$ is a lower triangular matrix
satisfying
$$
t_{ij}=
\begin{cases}
0,&\text{if } i < j;\\
1,&\text{if } i=j;\\
a_{i-j}(i-1),&\text{if } 2\leq i\leq n, 1\leq j\leq i-1.
\end{cases}
$$
Since $\mathbf{T}_{n}$ and $\mathbf{D}_{n}$ may contain too many parameters as compared with $n$,
we are led to consider a banded Cholesky decomposition of
$\Sigma_{n}^{-1}$,
\begin{equation}
\Sigma_{n}^{-1}(q)=\mathbf{T}_{n}'(q)\mathbf{D}_{n}^{-1}(q)\mathbf{T}_{n}(q),
\end{equation}
where $1\leq q\ll n$ is referred to as the banding parameter,
\begin{equation*}
\mathbf{D}_{n}(q)=diag(\gamma_{0}, \sigma_{1}^{2},\cdots, \sigma_{q}^{2}, \cdots, \sigma_{q}^{2}),
\end{equation*}
and
$\mathbf{T}_{n}(q)=(t_{ij}(q))_{1\leq i,j\leq n}$ with
$$
t_{ij}(q)=
\begin{cases}
0,&\text{if } i < j \ \text{or }  \{ q+1 < i \leq n, 1 \leq j
\leq
i-q-1 \};\\
1,&\text{if } i=j;\\
a_{i-j}(i-1),&\text{if } 2\leq i\leq q, 1\leq j\leq i-1;\\
a_{i-j}(q),&\text{if }q+1 \leq i \leq n, i-q \leq j \leq i-1.
\end{cases}
$$

To estimate the banded Cholesky factors in (3.3), we first generate the least squares residuals
$\hat{\mathbf{e}}_{n}=(\hat{e}_{1}, \ldots, \hat{e}_{n})^{'}$ based on the approximating model $\sum_{j=1}^{d}\theta_{j}x_{t, j}$
for (2.1),
where
$d=d_{n}$ is allowed to grow to infinity with $n$ and
$\hat{\mathbf{e}}_{n}=(I-H_{d})Y_{n}$ with $H_{d}$ denoting the orthogonal projection matrix for the column space of
$X(d)=(x_{ij})_{1\leq i \leq n, 1\leq j \leq d}$.
Having obtained
$\hat{\mathbf{e}}_{n}$,
the $\gamma_{0}$ and $\sigma^{2}(k)$ in $\mathbf{D}_{n}(q)$
and the $a_{i}(k)$ in $\mathbf{T}_{n}(q)$
can be estimated by
$\hat{\gamma}_{0}$, $\hat{\sigma}^{2}(k)$,
and
$\hat{a}_{i}(k)$, respectively, where for $1\leq k \leq q$,
\begin{eqnarray*}
(\hat{a}_{1}(k),\ldots,\hat{a}_{k}(k))^{'}&=&\text{arg$\min$}_{(c_{1},\ldots,c_{k})^{'}\in
R^{k}}\sum_{t=q+1}^{n}(\hat{e}_{t}+
c_{1}\hat{e}_{t-1}+\cdots+c_{k}\hat{e}_{t-k})^2,\\
\hat{\gamma}_{0}&=&n^{-1}\sum_{t=1}^{n}\hat{e}_{t}^{2},\\
\hat{\sigma}_{k}^{2}&=&(n-q)^{-1}\sum_{t=q+1}^{n}(\hat{e}_{t}+\sum_{j=1}^{k}
\hat{a}_{j}(k)\hat{e}_{t-j})^{2}.
\end{eqnarray*}
Plugging these estimators
into
$\mathbf{D}_{n}(q)$ and
$\mathbf{T}_{n}(q)$, we obtain $\hat{\mathbf{D}}_{n}(q)$
and
$\hat{\mathbf{T}}_{n}(q)$, and hence an estimator of $\Sigma^{-1}_{n}$,
\begin{eqnarray}
\hat{\Sigma}_{n}^{-1}(q)=\hat{\mathbf{T}}_{n}'(q)\hat{\mathbf{D}}_{n}^{-1}(q)\hat{\mathbf{T}}_{n}(q).
\end{eqnarray}
Note that
we have suppressed the dependence of $\hat{\Sigma}_{n}^{-1}(q)$ on $d$ in order to
simplify notation.
The next theorem provides a rate of convergence of $\hat{\Sigma}_{n}^{-1}(q)$ to $\Sigma^{-1}_{n}$
when $q=q_{n}$ and $d=d_{n}$ grow to infinity with $n$ at suitable rates.
\begin{thm}
Assume {\rm Assumptions 1 and 2}
with {\rm (2.13)} and {\rm (2.15)} replaced by
\begin{eqnarray}
\mathrm{E}(|\alpha_{t}|^{2r}|{\cal F}_{t-1})<C_{2r} <\infty \,\,{\rm a.s.},
\end{eqnarray}
where $r\geq 2$, and
\begin{eqnarray}
\sum_{j\geq 1}j|\beta_{j}|<\infty,
\end{eqnarray}
respectively. Also assume that
\begin{eqnarray}
\sup_{t\geq 1, j \geq 1}
\mathrm{E}(|x_{tj}|^{2r})<\infty.
\end{eqnarray}
Suppose that $d_{n}$ and $q_{n}$ are chosen to satisfy:
\begin{eqnarray}
d_{n}\asymp n^{1/4},
\end{eqnarray}
\begin{eqnarray}
\max{\{d_{n}, q_{n}\}}\sum_{j \geq d_{n}}|\theta_{j}|=o(1),
\end{eqnarray}
and
\begin{eqnarray}
\frac{q^{2}_{n}d_{n}}{n}=o(1).
\end{eqnarray}
Then,
\begin{eqnarray}
\|\hat{\Sigma}_{n}^{-1}(q_{n})-\Sigma^{-1}_{n}\|=O_{p}\left(\frac{q^{1+ r^{-1}}_{n}}{n^{1/2}}+
\sqrt{\sum_{j\geq q_{n}+1}|a_{j}|\sum_{j\geq q_{n}+1}j|a_{j}|}\right).
\end{eqnarray}
\end{thm}

\vspace{0.3cm}
\noindent
{\bf Remark 11.}
In the simpler situation where $\hat{\mathbf{e}}_{n}=Y_{n}=\mathbf{e}_{n}$, namely,
$\mu_{t}=0$ for all $t$, Wu and Pourahmadi (2009) proposed a banded covariance matrix estimator
$\breve{\Sigma}_{n, l}=(\hat{\gamma}_{i-j}I_{|i-j|\leq l})_{1\leq i, j\leq n}$
of $\Sigma_{n}$, where $\hat{\gamma}_{k}=n^{-1}\sum_{i=1}^{n-|k|}e_{i}e_{i+|k|}$
is the $k$th lag sample ACF of $\{e_{t}\}$ and
$l$ is also called the banding parameter.
When $l=l_{n}=o(n^{1/2})$ and (3.5) holds with $r=2$,
their Theorems 2 and 3 imply that
$\breve{\Sigma}_{n, l_{n}}$ is positive definite with probability approaching one,
\begin{eqnarray}
\|\breve{\Sigma}_{n, l_{n}}-\Sigma_{n}\|=O_{p}\left(\frac{l_{n}}{n^{1/2}}+
\sum_{j\geq q_{n}}|\gamma_{j}|\right),
\end{eqnarray}
and
\begin{eqnarray}
\|\breve{\Sigma}^{-1}_{n, l_{n}}-\Sigma^{-1}_{n}\|=O_{p}\left(\frac{l_{n}}{n^{1/2}}+
\sum_{j\geq q_{n}}|\gamma_{j}|\right).
\end{eqnarray}
McMurry and Politis (2010)
generalized (3.12) and (3.13) to tapered covariance matrix estimators.
Ing, Chiou and Guo (2013)
considered estimating $\Sigma^{-1}_{n}$ through the
banded Cholesky decomposition approach
in situations where
$\hat{\mathbf{e}}_{n}$ is obtained by
a {\it correctly specified} regression model.
They established the consistency of the proposed estimator
under spectral norm, even when $\{e_{t}\}$
is a long-memory time series.
However, since this section allows the regression model to be misspecified,
all the aforementioned results are not directly applicable here.
\\
\\
{\bf Remark 12.}
The second term on the right-hand side of (3.11) is mainly contributed by
the approximation error
$\|\Sigma_{n}^{-1}(q_{n})-\Sigma_{n}^{-1}\|$, whereas the first one is mainly due to
the sampling variability $\|\hat{\Sigma}_{n}^{-1}(q_{n})-\Sigma_{n}^{-1}(q_{n})\|$,
which is in turn dominated by
$\|\hat{\mathbf{T}}_{n}(q_{n})-\mathbf{T}_{n}(q_{n})\|$, as shown in the proof of Theorem 3.
Similarly,
the first and second terms on the right-hand side of (3.12)
are contributed by
$\|\breve{\Sigma}_{n, l_{n}}-\Sigma_{n, l_{n}}\|$
and $\|\Sigma_{n, l_{n}}-\Sigma_{n}\|$, respectively.
Here, $\Sigma_{n, l_{n}}=(\gamma_{i-j}I_{|i-j|\leq l_{n}})_{1\leq i, j\leq n}$
is the population version of $\breve{\Sigma}_{n, l_{n}}$.
However, unlike $\breve{\Sigma}_{n, l_{n}}-\Sigma_{n, l_{n}}$,
$\hat{\mathbf{T}}_{n}(q_{n})-\mathbf{T}_{n}(q_{n})$ is not a Toeplitz matrix.
Hence our upper bound for $\|\hat{\mathbf{T}}_{n}(q_{n})-\mathbf{T}_{n}(q_{n})\|$
is derived from complicated maximal probability inequalities, such as (A.42) and (A.49),
which also lead to an additional exponent $r^{-1}$ in the first term on the right-hand side of (3.11).
\\
\\
{\bf Remark 13.}
The technical assumptions (3.8)-(3.10) essentially say that
the dimension, $d_{n}$, of the working regression model shouldn't be too large or too small. They ensure that
the sampling variability and the approximation error
introduced by this model are completely absorbed into
the first or second term
on the right-hand side of (3.11),
which depend only on the working AR model
used in the Cholesky decomposition.
As shown in the next section,
this feature can substantially reduce the burden
of verifying (2.24) and (2.25).

\def \theequation{4.\arabic{equation}}
\setcounter{equation}{0}

\section{Asymptotic efficiency of the FAMMA method with $\hat{\Sigma}^{-1}_{n}=\hat{\Sigma}^{-1}_{n}(q_{n})$.}
In this section, we shall establish the asymptotic efficiency of
$\hat{C}^{*}_{n}(\boldsymbol{w})$
with $\hat{\Sigma}^{-1}_{n}=\hat{\Sigma}^{-1}_{n}(q_{n})$, denoted by $\hat{C}^{*}_{n,q_{n}}(\boldsymbol{w})$,
when the AR coefficients of $\{e_{t}\}$
satisfy
\begin{eqnarray}
\sqrt{\sum_{j\geq q}|a_{j}|\sum_{j\geq q}j|a_{j}|} \leq C_{1} {\rm exp}(-\nu q),
\end{eqnarray}
or
\begin{eqnarray}
\sqrt{\sum_{j\geq q}|a_{j}|\sum_{j\geq q}j|a_{j}|} \leq C_{2} q^{-\nu},
\end{eqnarray}
for all $q\geq 1$ and some positive constants $C_{1}, C_{2}$ and $\nu$.
We call (4.1) the exponential decay case, which is fulfilled by any causal and invertible
ARMA($p, q$) model with $0\leq p,q<\infty$.
On the other hand, (4.2)
is referred to as the algebraic decay case, which is commonly discussed in the context of model selection for time series;
see Shibata (1981) and Ing and Wei (2003, 2005).

We first choose suitable $q_{n}$ for $\hat{\Sigma}^{-1}_{n}(q_{n})$
to ensure that the bound in (3.11) possesses the optimal rate.
When (4.1) is assumed, it is not difficult to see that
the optimal rate of (3.11) is
$O_{p}((\log n)^{1+r^{-1}}/n^{1/2})$,
which is achieved by
\begin{eqnarray}
q_{n}= c_{4} \log n,
\end{eqnarray}
for some sufficiently large constant $c_{4}$.
Therefore, (2.24) holds
with
\begin{eqnarray}
\hat{\Sigma}^{-1}_{n}=\hat{\Sigma}^{-1}_{n}(c_{4} \log n) \,\,\mbox{and} \,\,b_{n}=(\log n)^{1+r^{-1}}.
\end{eqnarray}
When (4.2) is true,
by letting
\begin{eqnarray}
q_{n}= \lfloor n^{1/\{2(1+r^{-1}+\nu)\}}\rfloor,
\end{eqnarray}
where $\lfloor a \rfloor$ denotes the largest integer $\leq a$,
we get the optimal rate
of (3.11), $O_{p}(n^{-\nu/\{2(1+\nu+r^{-1})\}})$,
yielding that (2.24) holds
with
\begin{eqnarray}
\hat{\Sigma}^{-1}_{n}=\hat{\Sigma}^{-1}_{n}(\lfloor
n^{1/\{2(1+r^{-1}+\nu)\}}\rfloor) \,\,\mbox{and} \,\,b_{n}=n^{\frac{1+r^{-1}}{2(1+\nu+r^{-1})}}.
\end{eqnarray}

We are ready to establish the
asymptotic efficiency of
$\hat{C}^{*}_{n,q_{n}}(\boldsymbol{w})$ under (4.1).

\begin{col}
Assume {\rm Assumptions 1 and 2}, {\rm (4.1)} and {\rm (3.7)} with $2r$ replaced by $S$,
noting that $S$ is defined in {\rm (M3)} of {\rm Assumption 1}.
Suppose that $d_{n}$ and $q_{n}$ obey {\rm (3.8)} and {\rm (4.3)}, respectively.
Moreover, assume
\begin{eqnarray}
d_{n}\sum_{j \geq d_{n}}|\theta_{j}|=o(1),
\end{eqnarray}
and for some $2N/S<\theta<1/2$,
\begin{eqnarray}
\frac{(\log n)^{1+(2/S)}}{k^{*^{(1/2)-\theta}}_{n}}=0 \,\, {\rm a.s.}
\end{eqnarray}
Then, {\rm (2.26)} holds with $\hat{\boldsymbol{w}}_{n}=\arg\inf_{\boldsymbol{w} \in {\cal H}_{N}}\hat{C}^{*}_{n,q_{n}}(\boldsymbol{w})$.
\end{col}

Corollary 1 follows directly from Theorems 2 and 3 and (4.4) with $r^{-1}$ replaced by $2/S$. Its proof is thus omitted.
Condition (4.8) is easily satisfied when $D_{n}(m)$ follows (2.21).
To see this, note that
(2.21) implies $k^{*}_{n}=c_{5}n^{1/(1+a)}$ for some $c_{5}>0$.
Therefore, (4.8) holds for any $2N/S<\theta<1/2$.
On the other hand, it is not difficult to show that (4.8) is violated when
$D(m)=n{\rm exp}(-c_{6}m)+m$ for some $c_{6}>0$, which leads to a much smaller $k^{*}_{n}=c_{7} \log n$
for some $c_{7}>0$.

To establish the
asymptotic efficiency of
$\hat{C}^{*}_{n,q_{n}}(\boldsymbol{w})$ under (4.2) with $\nu$ unknown,
we need to assume that $\nu$ has a known lower limit $\nu_{0} \geq 1/3$.

\begin{col}
Assume {\rm Assumptions 1 and 2}, {\rm (4.2)} and {\rm (3.7)} with $2r$ replaced by $S$.
Suppose that $d_{n}$ obeys {\rm (3.8)} and $q_{n}$ satisfies {\rm (4.5)}
with $r^{-1}$ replaced by $2/S$ and $\nu$ by $\nu_{0}$.
Moreover, assume {\rm (3.9)}
and for some $2N/S<\theta<1/2$,
\begin{eqnarray}
\frac{n^{\frac{1+(2/S)}{2[1+\nu_{0}+(2/S)]}}}{k^{*^{(1/2)-\theta}}_{n}}=0 \,\, {\rm a.s.}
\end{eqnarray}
Then, the conclusion of Corollary 1 follows.
\end{col}
Corollary 2 can be proved using Theorems 2 and 3 and (4.6)
with $r^{-1}$ and $\nu$ replaced by $2/S$ and $\nu_{0}$, respectively.
We again omit the details.
Before closing this section, we provide
a sufficient condition for (4.9) in situations where
$D_{n}(m)$ obeys (2.21). We assume that
(2.13) in Assumption 1 holds for any $0<S<\infty$ in order to simplify exposition.
Elementary calculations show that (4.9) follows from
$\nu_{0}>a$. However, since $a$ is in general unknown,
our simple and practical guidance for verifying (4.9)
is to check whether $\nu_{0}>\bar{a}$, where $\bar{a}$ is a known upper bound for $a$.

\section{Concluding remarks}

This paper provides guidance for the model averaging implementation
in regression models with time series errors.
Driven by the efficiency improvement, our goal is to choose the optimal weight vector that averages across FGLS
estimators obtained from a set of approximating models of the true regression function.
We propose the FAMMA as the weight selection
criterion and show its asymptotic optimality in the sense of (2.26).
To the best of our knowledge, it is the first time that the
FGLS-based criterion is proved to have this type of property
in the presence of time-dependent errors.


On the other hand, our asymptotic optimality,
implicitly involving the search for the averaging estimator
whose loss (or conditional risk) has the best constant in addition to
the best rate, is typically not achievable when the number of candidate models is large and
the models are not necessarily nested.
While Wan, Zhang and Zou (2010) and Zhang, Wan and Zou (2013)
proved the asymptotic efficiency of their averaging estimators without assuming nested candidate models,
a stringent condition on the number of models, e.g., (2.20),
is placed as the tradeoff.
Furthermore, on top of their positive report, no clear guideline for
the optimal averaging across arbitrary combinations of regressors was offered.
In fact, in this more challenging situation,
pursuing the minimax optimal rate
appears to be more relevant than the asymptotic efficiency.
The theoretical results developed in Wang et al. (2014)
and in Sections 2 and 3 provide useful tools for deriving
the minimax optimal rate under model (2.1).
Moreover, motivated by Ing and Lai (2011), we conjecture that
when
the variables are preordered by
the orthogonal greedy algorithm (OGA) (see, e.g., Temlyakov (2000) and Ing and Lai (2011)),
this rate is achievable by
FAMMA with 2 replaced by
a factor directly proportional to
the natural logarithm of the number of candidate models.
We leave investigations along this research direction to future work.



\def \theequation{4.\arabic{equation}}
\setcounter{equation}{0}

\vspace{0.5cm}

\def \theequation{A.\arabic{equation}}
\setcounter{equation}{0}

\centerline{\bf\large APPENDIX}

\vspace{0.3cm}
{\it Proof of Lemma 1.}
Note first that
$R^{*}_{n}(\boldsymbol{w})=\mathrm{E}_{\boldsymbol{x}}(L^{*}_{n}(\boldsymbol{w}))=
\mathrm{E}_{\boldsymbol{x}}(\textbf{e}^{'}_{n}P^{*^{'}}(\boldsymbol{w})\Sigma^{-1}_{n}P^{*}(\boldsymbol{w})\textbf{e}_{n})+
\boldsymbol{\mu}^{'}_{n}(I-P^{*}(\boldsymbol{w}))^{'}\Sigma^{-1}_{n}(I-P^{*}(\boldsymbol{w}))\boldsymbol{\mu}_{n}$.
Since
\begin{eqnarray}
\Sigma^{-1/2}_{n}P^{*}(\boldsymbol{w})=\sum_{m=1}^{M}w_{m}P_{m}\Sigma^{-1/2}_{n},
\end{eqnarray}
it follows that
\begin{eqnarray*}
&&\mathrm{E}_{\boldsymbol{x}}(\textbf{e}^{'}_{n}P^{*^{'}}(\boldsymbol{w})\Sigma^{-1}_{n}P^{*}(\boldsymbol{w})\textbf{e}_{n})\\
&=&
\mathrm{E}_{\boldsymbol{x}}(\sum_{m=1}^{M}\sum_{l=1}^{M}w_{m}w_{l}\textbf{e}^{'}_{n}\Sigma^{-1/2}_{n}P_{l}P_{m}\Sigma^{-1/2}_{n}
\textbf{e}_{n})\\
&=&\sum_{m=1}^{M}\sum_{l=1}^{M}w_{l}w_{m}\min\{k_{m}, k_{l}\}.
\end{eqnarray*}
Similarly, $\boldsymbol{\mu}^{'}_{n}(I-P^{*}(\boldsymbol{w}))^{'}\Sigma^{-1}_{n}(I-P^{*}(\boldsymbol{w}))\boldsymbol{\mu}_{n}=
\sum_{m=1}^{M}\sum_{l=1}^{M}w_{l}w_{m}\boldsymbol{\mu}^{'}_{n}\Sigma^{-1/2}_{n}(I-P_{\max\{m, l\}})\Sigma^{-1/2}_{n}\boldsymbol{\mu}_{n}$.
Consequently, the desired conclusion (2.2) follows.



\vspace{0.3cm}
\noindent
{\it Proof of Theorem 1.}
Define  $\boldsymbol{w}^{*}_{n} = \arg
\displaystyle{\min_{\boldsymbol{w} \in \mathcal{H}_{N}}} L_{n}^{*} ( \boldsymbol{w} )$,
$(\tilde{w}_{n, 1}, \ldots, \tilde{w}_{n, M})^{'}=\tilde{\boldsymbol{w}}_{n}$,
and
$(w^{*}_{n, 1}, \ldots, w^{*}_{n, M})^{'}=\boldsymbol{w}^{*}_{n}$.
By noticing
\begin{eqnarray*}
&&C^{*}_{n}(\boldsymbol{w})-L^{*}_{n}(\boldsymbol{w})=
\textbf{e}^{'}_{n}\Sigma^{-1}_{n}\textbf{e}_{n}+
2\textbf{e}^{'}_{n}\Sigma^{-1}_{n}(I-P^{*}(\boldsymbol{w}))\boldsymbol{\mu}_{n} \nonumber \\
&-&2\{\textbf{e}^{'}_{n}\Sigma^{-1}_{n}P^{*}(\boldsymbol{w})\textbf{e}_{n}-\sum_{m=1}^{M}w_{m}k_{m}\},
\end{eqnarray*}
we get
\begin{eqnarray}
0 &\geq&  \left\{ C_{n}^{*} ( \tilde{\boldsymbol{w}}_{n} ) -  C_{n}^{*} ( \boldsymbol{w}^{*}_{n} ) \right\} =
 L_{n}^{*} ( \tilde{\boldsymbol{w}}_{n} ) -  L_{n}^{*} ( \boldsymbol{w}^{*}_{n} )  +
2 \mathbf{e}'_{n} \Sigma_{n}^{-1} \left( I - P^{*} ( \tilde{\boldsymbol{w}}_{n} ) \right) \boldsymbol{\mu}_{n} \nonumber \\
&-& 2 \left\{\mathbf{e}'_{n} \Sigma_{n}^{-1}
P^{*} ( \tilde{\boldsymbol{w}}_{n} ) \mathbf{e}_{n} -
\sum_{m=1}^{M} \tilde{w}_{n, m} k_{m}   \right\} - 2 \mathbf{e}'_{n} \Sigma_{n}^{-1}
\left( I - P^{*} ( \boldsymbol{w}^{*}_{n} ) \right) \boldsymbol{\mu}_{n} \nonumber \\
&+& 2 \left\{ \mathbf{e}'_{n} \Sigma_{n}^{-1} P^{*} ( \boldsymbol{w}^{*}_{n} )
\mathbf{e}_{n} - \sum_{m=1}^{M} w^{*}_{n, m} k_{m} \right\} \nonumber \\
&=& L_{n}^{*} ( \tilde{\boldsymbol{w}}_{n} ) -  L_{n}^{*} ( \boldsymbol{w}^{*}_{n} )
+ 2 A_{n} ( \tilde{\boldsymbol{w}}_{n} ) - 2 B_{n} ( \tilde{\boldsymbol{w}}_{n} ) - 2
A_{n} ( \boldsymbol{w}^{*}_{n} ) + 2 B_{n} ( \boldsymbol{w}^{*}_{n} ),
\end{eqnarray}
where $A_{n} ( \boldsymbol{w} ) = \mathbf{e}'_{n} \Sigma_{n}^{-1} \left( I - P^{*} ( \boldsymbol{w} )
\right) \boldsymbol{\mu}_{n}$ and $B_{n} ( \boldsymbol{w} ) = \mathbf{e}'_{n} \Sigma_{n}^{-1} P^{*} ( \boldsymbol{w} )
\mathbf{e}_{n} - \sum_{m=1}^{M} w_{m} k_{m}$.
In view of (A.2) and $L_{n}^{*} ( \tilde{\boldsymbol{w}}_{n} ) \geq  L_{n}^{*} ( \boldsymbol{w}^{*}_{n} )$,
it suffices for (2.11) to show that
\begin{eqnarray}
\displaystyle{\sup_{\boldsymbol{w} \in \mathcal{H}_{N}}} \left|
\displaystyle{\frac{A_{n} ( \boldsymbol{w} )}{ R^{*}_{n} ( \boldsymbol{w} )}} \right| = o_{p} (1),
\end{eqnarray}
\begin{eqnarray}
\displaystyle{\sup_{\boldsymbol{w} \in \mathcal{H}_{N}}} \left|
\displaystyle{\frac{B_{n} ( \boldsymbol{w} )}{R^{*}_{n} ( \boldsymbol{w} )}} \right| = o_{p} (1),
\end{eqnarray}
and
\begin{eqnarray}
\displaystyle{\sup_{\boldsymbol{w} \in \mathcal{H}_{N}}}
\left|\displaystyle{\frac{L^{*}_{n} ( \boldsymbol{w} )}{R^{*}_{n} ( \boldsymbol{w} )}}-1\right|=o_{p}(1),
\end{eqnarray}
where $\stackrel{p}{\longrightarrow}$ denotes convergence in probability.

To show (A.3), first note that
\begin{eqnarray*}
\mathcal{H}_{( l )} = \bigcup_{1 \leq j_{1} < j_{2} < \cdots < j_{l} \leq M} \mathcal{H}_{j_{1}, \cdots, j_{l}},
\end{eqnarray*}
where for $1 \leq j_{1} < \cdots < j_{l} \leq M$, $\mathcal{H}_{j_{1}, \cdots, j_{l}} = \{ \boldsymbol{w} :
\boldsymbol{w} \in \mathcal{H}_{( l )} \, {\rm and} \, \omega_{j_{i}} \neq 0, 1 \leq i \leq l \}$. Hence for any
$\varepsilon > 0$,
\begin{eqnarray}
& & P_{\boldsymbol{x}} \left( \displaystyle{\sup_{\boldsymbol{w} \in \mathcal{H}_{N}}} \left|
\displaystyle{\frac{A_{n} ( \boldsymbol{w} )}{ R^{*}_{n} ( \boldsymbol{w} )}} \right| > \varepsilon \right) \nonumber \\
&\leq& \sum_{l=1}^{N} \sum_{j_{l}=l}^{M} \cdots \sum_{j_{1}=1}^{j_{2}-1}
P_{\boldsymbol{x}} \left( \displaystyle{\sup_{\boldsymbol{w} \in \mathcal{H}_{j_{1}, \cdots, j_{l}}}} \left|
\displaystyle{\frac{A_{n} ( \boldsymbol{w} )}{ R^{*}_{n} ( \boldsymbol{w} )}} \right| > \varepsilon \right) \nonumber \\
&\leq& \sum_{l=1}^{N} \sum_{j_{l}=l}^{M} \cdots \sum_{j_{1}=1}^{j_{2}-1} P_{\boldsymbol{x}} \left(
\displaystyle{\frac{\displaystyle_{\sum_{m \in \{ j_{1}, \cdots, j_{l} \}}} \left| \boldsymbol{\mu}'_{n} \Sigma_{n}^{-1/2}
( I - P_{m} ) \Sigma_{n}^{-1/2} \mathbf{e}_{n} \right|}
{\underline{\delta}^{2} \displaystyle{\max_{m \in \{ j_{1}, \cdots, j_{l} \}}} D_{n} ( m )}} > \varepsilon \right)\equiv
 \sum_{l=1}^{N} Q_{l},
\end{eqnarray}
where $P_{\boldsymbol{x}} (\cdot)=P(\cdot|\boldsymbol{x}_{1}, \ldots, \boldsymbol{x}_{n})$
and the second inequality follows from
\begin{eqnarray}
\displaystyle{\inf_{\boldsymbol{w} \in \mathcal{H}_{j_{1}, \cdots, j_{l}}}}  R_{n}^{*} ( \boldsymbol{w} ) \geq
\underline{\delta}^{2} \displaystyle{\max_{m \in \{ j_{1}, \cdots, j_{l} \}}} D_{n} ( m ),
\end{eqnarray}
which is ensured by Lemma 1 and the definition of ${\cal H}_{N}$.
Let $S_{1}=S/2$. Then, by Chebshev's inequality, (M3), and Lemma 2 of Wei (1987), it holds that
\begin{eqnarray*}
& & P_{\boldsymbol{x}} \left( \displaystyle{\frac{\displaystyle_{\sum_{m \in \{ j_{1}, \cdots, j_{l} \}}} \left| \boldsymbol{\mu}'_{n} \Sigma_{n}^{-1/2}
( I - P_{m} ) \Sigma_{n}^{-1/2} \mathbf{e}_{n} \right|}
{\underline{\delta}^{2} \displaystyle{\max_{m \in \{ j_{1}, \cdots, j_{l} \}}} D_{n} ( m )}} > \varepsilon \right) \nonumber \\ &\leq&
C \sum_{m \in \{ j_{1}, \cdots, j_{l} \}} \left( \displaystyle
\frac{\mathrm{E}_{\boldsymbol{x}}\{ \boldsymbol{\mu}'_{n} \Sigma_{n}^{-1/2}
( I - P_{m} ) \Sigma_{n}^{-1/2} \mathbf{e}_{n} \}^{S_{1}}}
{\{\displaystyle{\max_{m \in \{ j_{1}, \cdots, j_{l} \}}} D_{n} ( m )\}^{S_{1}}}     \right) \nonumber \\
&\leq& C \sum_{m \in \{ j_{1}, \cdots, j_{l} \}} \displaystyle{\frac{\left( \boldsymbol{\mu}'_{n} \Sigma_{n}^{-1/2}
( I - P_{m} ) \Sigma_{n}^{-1/2} \boldsymbol{\mu}_{n} \right)^{S_{1}/2}}{\left( \displaystyle{\max_{m \in \{ j_{1}, \cdots, j_{l} \}}}
D_{n} ( m ) \right)^{S_{1}}}} \nonumber \\ &\leq& C \sum_{m \in \{ j_{1}, \cdots, j_{l} \}} \displaystyle{\frac{1}
{\left( \displaystyle{\max_{m \in \{ j_{1}, \cdots, j_{l} \}}} D_{n} ( m ) \right)^{S_{1}/2}}}
\leq C \displaystyle{\frac{l}{\left( D_{n} ( j_{l} ) \right)^{S_{1}/2}}},
\end{eqnarray*}
where here and hereafter $C$ denotes a generic positive constant whose value is independent of
$n$ and may vary at different occurrences.
Therefore, for each $1 \leq l \leq N$,
\begin{eqnarray*}
Q_{l} &\leq& C \left\{ \sum_{j_{l}=l}^{k_{n}^{*}} \cdots \sum_{j_{1}=1}^{j_{2}-1} \displaystyle{\frac{1}{\left( k_{n}^{*} \right)^{S_{1}/2}}}
+ \sum_{j_{l}=k_{n}^{*}+1}^{M} \cdots \sum_{j_{1}=1}^{j_{2}-1} \displaystyle{\frac{1}{\left( D_{n} ( j_{l} ) \right)^{S_{1}/2}}} \right\}
\\ &\leq& C \left\{ k_{n}^{{*}^{-(S_{1}/2-l)}} +
\sum_{j_{l}=k_{n}^{*}+1}^{\infty} \displaystyle{\frac{j_{l}^{l-1}}{j_{l}^{S_{1}/2}}} \right\},
\end{eqnarray*}
which converges to 0 a.s. in view of (2.16).
As a result,
\begin{eqnarray*}
P_{\boldsymbol{x}} \left( \displaystyle{\sup_{\boldsymbol{w} \in \mathcal{H}_{N}}} \left|
\displaystyle{\frac{A_{n} ( \boldsymbol{w} )}{ R^{*}_{n} ( \boldsymbol{w} )}} \right| > \varepsilon \right)
\rightarrow 0, \, {\rm a.s.}
\end{eqnarray*}
This and the dominated convergence theorem together imply (A.3).

Similarly,
\begin{eqnarray*}
P_{\boldsymbol{x}} \left( \displaystyle{\sup_{\boldsymbol{w} \in \mathcal{H}_{N}}} \left|
\displaystyle{\frac{B_{n} ( \boldsymbol{w} )}{ R^{*}_{n} ( \boldsymbol{w} )}} \right| > \varepsilon \right)
\leq C\sum_{l=1}^{N} E_{l},
\end{eqnarray*}
where
\begin{eqnarray*}
E_{l} = \sum_{j_{l}=l}^{M} \cdots \sum_{j_{1}=1}^{j_{2}-1} P_{\boldsymbol{x}} \left(
\displaystyle{\frac{\displaystyle_{\sum_{m \in \{ j_{1}, \cdots, j_{l} \}}} \left| \mathbf{e}'_{n} \Sigma_{n}^{-1/2}
P_{m} \Sigma_{n}^{-1/2} \mathbf{e}_{n} - k_{m} \right|}
{\underline{\delta}^{2} \displaystyle{\max_{m \in \{ j_{1}, \cdots, j_{l} \}}} D_{n} ( m )}} > \varepsilon \right).
\end{eqnarray*}
By (M3) and the first moment bound theorem of Findley and Wei (1993), it follows that
\begin{eqnarray*}
&&P_{\boldsymbol{x}} \left( \displaystyle{\frac{\displaystyle_{\sum_{m \in \{ j_{1}, \cdots, j_{l} \}}} \left| \mathbf{e}'_{n} \Sigma_{n}^{-1/2}
P_{m} \Sigma_{n}^{-1/2} \mathbf{e}_{n} - k_{m} \right|}
{\underline{\delta}^{2} \displaystyle{\max_{m \in \{ j_{1}, \cdots, j_{l} \}}} D_{n} ( m )}} > \varepsilon \right) \\
&\leq&
C \sum_{m \in \{ j_{1}, \cdots, j_{l} \}} \displaystyle{\frac{k_{m}^{S_{1}/2}}{\left(
\displaystyle{\max_{m \in \{ j_{1}, \cdots, j_{l} \}}} D_{n} ( m ) \right)^{S_{1}}}} \leq
\displaystyle{\frac{C l}{\left( D_{n} ( j_{l} ) \right)^{S_{1}/2}}}.
\end{eqnarray*}
Therefore, (A.4) follows immediately from an argument similar to that used to prove (A.3).
The proof of (A.5) is similar to those of (A.3) and (A.4). The details are omitted.

\vspace{0.5cm}
\noindent
{\it Proof of Theorem 2}.
Define $\hat{L}_{n}^{*} ( \boldsymbol{w} ) = \left( \hat{\boldsymbol{\mu}}_{n}^{*} ( \boldsymbol{w} ) - \boldsymbol{\mu}_{n} \right)'
\hat{\Sigma}_{n}^{-1} \left( \hat{\boldsymbol{\mu}}_{n}^{*} ( \boldsymbol{w} ) - \boldsymbol{\mu}_{n} \right)$.
Then, it follows that
\begin{eqnarray*}
\hat{C}_{n}^{*} ( \boldsymbol{w} ) &=& \left( \mathbf{e}_{n} - \left( \hat{\boldsymbol{\mu}}_{n}^{*} ( \boldsymbol{w} )
- \boldsymbol{\mu}_{n} \right) \right)' \hat{\Sigma}_{n}^{-1} \left( \mathbf{e}_{n} - \left( \hat{\boldsymbol{\mu}}_{n}^{*} ( \boldsymbol{w} )
- \boldsymbol{\mu}_{n} \right) \right) + 2 \sum_{m=1}^{M} w_{m} k_{m} \\
&=& \mathbf{e}'_{n} \hat{\Sigma}_{n}^{-1} \mathbf{e}_{n} +  \hat{L}_{n}^{*} ( \boldsymbol{w} ) - 2
\left( \hat{\boldsymbol{\mu}}_{n}^{*} ( \boldsymbol{w} ) - \boldsymbol{\mu}_{n} \right)' \hat{\Sigma}_{n}^{-1} \mathbf{e}_{n} +
2 \sum_{m=1}^{M} w_{m} k_{m} \\
&=& \mathbf{e}'_{n} \hat{\Sigma}_{n}^{-1} \mathbf{e}_{n} +
2 \boldsymbol{\mu}'_{n} \left( I - \hat{P}^{*} ( \boldsymbol{w} ) \right)'
\hat{\Sigma}_{n}^{-1} \mathbf{e}_{n} +  \hat{L}_{n}^{*} ( \boldsymbol{w} ) - 2 \left\{
\mathbf{e}'_{n} \hat{P}^{*'} ( \boldsymbol{w} ) \hat{\Sigma}_{n}^{-1} \mathbf{e}_{n} - \sum_{m=1}^{M} w_{m} k_{m} \right\},
\end{eqnarray*}
and hence
\begin{eqnarray*}
0 &\geq&  \hat{C}_{n}^{*} ( \hat{\boldsymbol{w}}_{n} ) - \hat{C}_{n}^{*} ( \boldsymbol{w}_{n}^{F} )   \\
&=&  \hat{L}_{n}^{*} ( \hat{\boldsymbol{w}}_{n} ) - \hat{L}_{n}^{*} ( \boldsymbol{w}_{n}^{F} )  +
2 \hat{A}_{n} ( \hat{\boldsymbol{w}}_{n} ) - 2 \hat{A}_{n} ( \boldsymbol{w}_{n}^{F} ) -
2 \hat{B}_{n} ( \hat{\boldsymbol{w}}_{n} ) + 2 \hat{B}_{n} ( \boldsymbol{w}_{n}^{F} ) \\
&=&  ( \hat{L}_{n}^{*} ( \hat{\boldsymbol{w}}_{n} ) - L^{F}_{n} ( \hat{\boldsymbol{w}}_{n} ) ) -
( \hat{L}_{n}^{*} ( \boldsymbol{w}_{n}^{F} ) - L_{n}^{F} ( \boldsymbol{w}_{n}^{F} ) ) +
2 \hat{A}_{n} ( \hat{\boldsymbol{w}}_{n} ) - 2 \hat{A}_{n} ( \boldsymbol{w}_{n}^{F} ) \\
&-& 2 \hat{B}_{n} ( \hat{\boldsymbol{w}}_{n} ) + 2 \hat{B}_{n} ( \boldsymbol{w}_{n}^{F} ) +
 ( L_{n}^{F} ( \hat{\boldsymbol{w}}_{n} ) - L_{n}^{F} ( \boldsymbol{w}_{n}^{F} ) ),
\end{eqnarray*}
where
$\boldsymbol{w}^{F}_{n}=\arg\min_{\boldsymbol{w} \in {\cal H}_{N}}L^{F}_{n}(\boldsymbol{w})$,
$\hat{A}_{n} ( \boldsymbol{w} )=\boldsymbol{\mu}'_{n} \left( I - \hat{P}^{*} ( \boldsymbol{w} ) \right)'
\hat{\Sigma}_{n}^{-1} \mathbf{e}_{n}$ and $\hat{B}_{n} ( \boldsymbol{w} )=
\mathbf{e}'_{n} \hat{P}^{*'} ( \boldsymbol{w} ) \hat{\Sigma}_{n}^{-1} \mathbf{e}_{n} - \sum_{m=1}^{M} w_{m} k_{m}$.
Since
$L_{n}^{F} ( \hat{\boldsymbol{w}}_{n} ) \geq L_{n}^{F} ( \boldsymbol{w}_{n}^{F} )$ and
(A.3)-(A.5) hold under the assumptions of Theorem 2, it suffices for (2.26) to show that
\begin{eqnarray}
\displaystyle{\sup_{\boldsymbol{w} \in \mathcal{H}_{N}}} \left|
\displaystyle{\frac{\hat{L}_{n}^{*} ( \boldsymbol{w} ) - L^{*}_{n} ( \boldsymbol{w} )}{R_{n}^{*} ( \boldsymbol{w} )}}
\right| = o_{p} ( 1 ),
\end{eqnarray}
\begin{eqnarray}
 \displaystyle{\sup_{\boldsymbol{w} \in \mathcal{H}_{N}}} \left|
\displaystyle{\frac{\hat{A}_{n} ( \boldsymbol{w} ) - A_{n} ( \boldsymbol{w} )}{ R_{n}^{*} ( \boldsymbol{w} )}}
\right| = o_{p} ( 1 ),
\end{eqnarray}
\begin{eqnarray}
\displaystyle{\sup_{\boldsymbol{w} \in \mathcal{H}_{N}}} \left|
\displaystyle{\frac{\hat{B}_{n} ( \boldsymbol{w} ) - B_{n} ( \boldsymbol{w} )}{ R_{n}^{*} ( \boldsymbol{w} )}}
\right| = o_{p} ( 1 ),
\end{eqnarray}
and
\begin{eqnarray}
\displaystyle{\sup_{\boldsymbol{w} \in \mathcal{H}_{N}}} \left|
\displaystyle{\frac{L_{n}^{F} ( \boldsymbol{w} ) - L^{*}_{n} ( \boldsymbol{w} )}{R_{n}^{*} ( \boldsymbol{w} )}}
\right| = o_{p} ( 1 ).
\end{eqnarray}

To prove (A.9), note first that
\begin{eqnarray}
&&\left| \hat{A}_{n} ( \boldsymbol{w} ) - A_{n} ( \boldsymbol{w} ) \right| \leq \left| \boldsymbol{\mu}_{n}'
( I - P^{*} ( \boldsymbol{w} ) )' ( \Sigma_{n}^{-1} - \hat{\Sigma}_{n}^{-1} ) \mathbf{e}_{n} \right| \nonumber \\
&+& \left| \boldsymbol{\mu}_{n}' ( \hat{P}^{*} ( \boldsymbol{w} ) - P^{*} ( \boldsymbol{w} ) )'
( \hat{\Sigma}_{n}^{-1} - \Sigma_{n}^{-1} ) \mathbf{e}_{n} \right| + \left| \boldsymbol{\mu}_{n}'
( \hat{P}^{*} ( \boldsymbol{w} ) - P^{*} ( \boldsymbol{w} ) )' \Sigma_{n}^{-1} \mathbf{e}_{n} \right| \nonumber \\
&\equiv& {\rm (1)} + {\rm (2)} + {\rm (3)}.
\end{eqnarray}
Assumption 2 implies
\begin{eqnarray}
\sup_{n \geq 1}\|\Sigma^{-1}_{n}\|<\infty \,\,\mbox{and} \,\,\sup_{n \geq 1}\|\Sigma_{n}\|<\infty,
\end{eqnarray}
which, together with (2.24), gives
\begin{eqnarray*}
&&(1) \leq C \| \Sigma_{n}^{-1} - \hat{\Sigma}_{n}^{-1}\| \|\mathbf{e}_{n}\|
 \| \Sigma_{n}^{-1/2}( I - P^{*} ( \boldsymbol{w} ) )\boldsymbol{\mu}_{n} \| \\
&=& O_{p} ( b_{n} ) R^{*^{1/2}}_{n}(\boldsymbol{w}),
\end{eqnarray*}
where the $O_{p} ( b_{n} )$ term is independent of $\boldsymbol{w}$.
In view of this, (A.7) and (2.25), one obtains
\begin{eqnarray}
&&
\displaystyle{\sup_{\boldsymbol{w} \in \mathcal{H}_{N}}}
\displaystyle{\frac{\left| \boldsymbol{\mu}_{n}' ( I - P^{*} ( \boldsymbol{w} ) )' ( \Sigma_{n}^{-1} - \hat{\Sigma}_{n}^{-1} )
\mathbf{e}_{n} \right|}{ R_{n}^{*} ( \boldsymbol{w} )}} \nonumber \\
&=&
\displaystyle{\max_{1\leq l \leq N}\max_{1\leq j_{1}<\cdots<j_{l}\leq M}\sup_{\boldsymbol{w} \in \mathcal{H}_{j_{1}, \cdots, j_{l}}}}
\displaystyle{\frac{\left| \boldsymbol{\mu}_{n}' ( I - P^{*} ( \boldsymbol{w} ) )' ( \Sigma_{n}^{-1} - \hat{\Sigma}_{n}^{-1} )
\mathbf{e}_{n} \right|}{ R_{n}^{*} ( \boldsymbol{w} )}} \nonumber \\
&=& O_{p} ( b_{n} ) \displaystyle{\frac{1}{k_{n}^{*^{1/2}}}}=o_{p}(1).
\end{eqnarray}
Let $A_{m}=X^{'}_{m}\Sigma^{-1}_{n}X_{m}$
and $\hat{A}_{m}=X^{'}_{m}\hat{\Sigma}^{-1}_{n}X_{m}$.
Then, straightforward calculations yield for any $\boldsymbol{w} \in \mathcal{H}_{j_{1}, \cdots, j_{l}}$
with $1\leq j_{1}<\cdots<j_{l}\leq M$ and $1\leq l \leq N$,
\begin{eqnarray}
{\rm (3)} \leq \sum_{m \in \{j_{1}, \ldots, j_{l}\}}
\left| \boldsymbol{\mu}_{n}' \left( \hat{\Sigma}_{n}^{-1} X_{m} ( \hat{A}_{m}^{-1} - A_{m}^{-1} )
X_{m}' + ( \hat{\Sigma}_{n}^{-1} - \Sigma^{-1}_{n} ) X_{m} A_{m}^{-1} X_{m}' \right) \Sigma_{n}^{-1} \mathbf{e}_{n} \right|.
\end{eqnarray}
It follows from (A.13) and (2.24) that
\begin{eqnarray*}
\sum_{m \in \{j_{1}, \ldots, j_{l}\}}
\left| \boldsymbol{\mu}_{n}' ( \hat{\Sigma}_{n}^{-1} - \Sigma^{-1} ) X_{m} A^{-1}_{m} X_{m}'
\Sigma_{n}^{-1} \mathbf{e}_{n} \right|
= O_{p} ( b_{n} ) \sum_{m \in \{j_{1}, \ldots, j_{l}\}}
\left\| P_{m} \Sigma_{n}^{-1/2}  \mathbf{e}_{n} \right\|.
\end{eqnarray*}
In addition, (A.7) and the first moment bound theorem of Findley and Wei (1993) imply that for $S_{1}=S/2>N/\theta$,
\begin{eqnarray*}
&&P_{\boldsymbol{x}} \left(
\displaystyle{\max_{1\leq l \leq N}\max_{1\leq j_{1}<\cdots<j_{l}\leq M}
\sup_{\boldsymbol{w} \in \mathcal{H}_{j_{1},\ldots,j_{l}}}} \displaystyle{\frac{\sum_{m \in \{j_{1},\ldots,j_{l}\}} \left\|
P_{m} \Sigma_{n}^{-1/2}  \mathbf{e}_{n} \right\|}{ R_{n}^{*} ( \boldsymbol{w} )}} > k_{n}^{*^{-1/2 + \theta}} \right)\\
&\leq& C \cdot k_{n}^{*^{S_{1}/2}} k_{n}^{*^{-S_{1} \theta}} \sum_{l=1}^{N} \sum_{j_{l}=l}^{M} \cdots \sum_{j_{1}=1}^{j_{2}-1}
\displaystyle{\frac{l}{D_{n}^{S_{1}/2} ( j_{l} )}} \\ &\leq& C \cdot k_{n}^{*^{S_{1}/2}} k_{n}^{*^{-S_{1} \theta}} k_{n}^{*^{-S_{1}/2+N}}.
\end{eqnarray*}
Combining the above two equations with (2.25) and the dominated convergence theorem, we get
\begin{eqnarray}
\displaystyle{\max_{1\leq l \leq N}\max_{1\leq j_{1}<\cdots<j_{l}\leq M}\sup_{\boldsymbol{w} \in \mathcal{H}_{j_{1}, \cdots, j_{l}}}}
\frac{\sum_{m \in \{ j_{1}, \cdots, j_{l} \}} \left|
\boldsymbol{\mu}_{n}' ( \hat{\Sigma}_{n}^{-1} - \Sigma^{-1} ) X_{m} A^{-1} X_{m}' \Sigma_{n}^{-1} \mathbf{e}_{n} \right|}
{R^{*}_{n}(\boldsymbol{w})} = o_{p} ( 1 ).
\end{eqnarray}
Some algebraic manipulations yield
\begin{eqnarray*}
\sum_{m \in \{ j_{1}, \cdots, j_{l} \}} \left| \boldsymbol{\mu}_{n}' \hat{\Sigma}_{n}^{-1} X_{m} ( \hat{A}_{m}^{-1} - A_{m}^{-1} )
X_{m}' \Sigma_{n}^{-1} \mathbf{e}_{n} \right| = O_{p}( b_{n}) \sum_{m \in \{ j_{1}, \cdots, j_{l} \}} \left\|
P_{m} \Sigma_{n}^{-1/2} \mathbf{e}_{n} \right\|.
\end{eqnarray*}
Therefore,
by an argument similar to that used to prove (A.16),
\begin{eqnarray}
\displaystyle{\max_{1\leq l \leq N}\max_{1\leq j_{1}<\cdots<j_{l}\leq M}\sup_{\boldsymbol{w} \in \mathcal{H}_{j_{1}, \cdots, j_{l}}}}
\displaystyle{\frac{\sum_{m \in \{ j_{1}, \cdots, j_{l} \}} \left| \boldsymbol{\mu}_{n}' \hat{\Sigma}_{n}^{-1} X_{m}
( \hat{A}_{m}^{-1} - A_{m}^{-1} ) X_{m}' \Sigma_{n}^{-1} \mathbf{e}_{n} \right|}{R_{n}^{*} ( \boldsymbol{w} )}} = o_{p} ( 1 ).
\end{eqnarray}
We conclude from (A.15), (A.16) and (A.17) that
\begin{eqnarray}
\displaystyle{\sup_{\boldsymbol{w} \in \mathcal{H}_{N}}} \displaystyle{\frac{\left| \boldsymbol{\mu}_{n}'
( \hat{P}^{*} ( \boldsymbol{w} ) - P^{*} ( \boldsymbol{w} ) )' \Sigma_{n}^{-1} \mathbf{e}_{n} \right|}
{R_{n}^{*} ( \boldsymbol{w} )}} =o_{p} ( 1 ).
\end{eqnarray}

Finally, straightforward calculations and (2.24) yield that
for any $\boldsymbol{w} \in \mathcal{H}_{j_{1}, \cdots, j_{l}}$
with $1\leq j_{1}<\cdots<j_{l}\leq M$ and $1\leq l \leq N$,
\begin{eqnarray*}
& & (2) \leq \displaystyle{\sum_{m \in \{ j_{1}, \cdots, j_{l} \}}} \left| \boldsymbol{\mu}_{n}' \hat{\Sigma}_{n}^{-1}
X_{m} ( \hat{A}_{m}^{-1} - A_{m}^{-1} ) X_{m}' ( \hat{\Sigma}_{n}^{-1} - \Sigma_{n}^{-1} ) \mathbf{e}_{n} \right| \\
&+& \displaystyle{\sum_{m \in \{ j_{1}, \cdots, j_{l} \}}} \left| \boldsymbol{\mu}_{n}'
( \hat{\Sigma}_{n}^{-1} - \Sigma_{n}^{-1} ) X_{m} ( \hat{A}_{m}^{-1} - A_{m}^{-1} ) X_{m}'
( \hat{\Sigma}_{n}^{-1} - \Sigma_{n}^{-1} ) \mathbf{e}_{n} \right| \\
&=& O_{p}(b^{2}_{n}).
\end{eqnarray*}
This and (2.25) imply
\begin{eqnarray}
\displaystyle{\sup_{\boldsymbol{w} \in \mathcal{H}_{N}}} \displaystyle{\frac{\left| \boldsymbol{\mu}_{n}'
( \hat{P}^{*} ( \boldsymbol{w} ) - P^{*} ( \boldsymbol{w} ) )' ( \hat{\Sigma}_{n}^{-1} - \Sigma_{n}^{-1} )
\mathbf{e}_{n} \right|}{ R_{n}^{*} ( \boldsymbol{w} )}} =O_{p}(b^{2}_{n}/k^{*}_{n})= o_{p} ( 1 ).
\end{eqnarray}
Now the desired conclusion (A.9) follows from (A.12), (A.14), (A.18) and (A.19).
The proofs of (A.8), (A.10) and (A.11) are similar to that of (A.9). The details are thus skipped.

\vspace{0.3cm}
Before proving Theorem 3, we need an auxiliary lemma.
\begin{lam}
Assume {\rm (2.12), (2.14)} and {\rm (3.6)}.
Then for any $1\leq q \leq n-1$,
\begin{eqnarray}
\|\Sigma^{-1}_{n}(q)-\Sigma^{-1}_{n}\|\leq C\sqrt{\sum_{j\geq q+1}|a_{j}|\sum_{j\geq q+1}j|a_{j}|},
\end{eqnarray}
where $a_{j}$'s are defined as in {\rm (3.1)}.
\end{lam}

{\it Proof.}
It follows from (2.12), (2.14), (3.6) and Theorem 3.8.4 of Brillinger (1975) that
\begin{eqnarray}
\sum_{j\geq 1}j|a_{j}|<\infty.
\end{eqnarray}
In view of (3.2) and (3.3), one has
\begin{eqnarray}
&&\|\Sigma_{n}^{-1}(q)-\Sigma^{-1}_{n}\|\leq
\|\mathbf{T}_{n}-\mathbf{T}_{n}(q)\|\|\mathbf{D}^{-1}_{n}\|\|\mathbf{T}_{n}\|+
\|\mathbf{T}_{n}(q)\|\|\mathbf{D}^{-1}_{n}-\mathbf{D}^{-1}_{n}(q)\|\|\mathbf{T}_{n}\| \nonumber \\
&+&
\|\mathbf{T}_{n}(q)\|\|\mathbf{D}^{-1}_{n}(q)\|\|\mathbf{T}_{n}-\mathbf{T}_{n}(q)\|.
\end{eqnarray}
It is easy to see that
\begin{eqnarray}
\|\mathbf{D}^{-1}_{n}(q)\| \leq C \,\,\mbox{and} \,\,\|\mathbf{D}^{-1}_{n}\| \leq C.
\end{eqnarray}
Moreover, by (A.13) and (A.23)
\begin{eqnarray}
\|\mathbf{T}_{n}\| \leq \|\mathbf{T}^{'}_{n}\mathbf{D}^{-1}_{n}\mathbf{T}_{n}\|\|\mathbf{D}_{n}\|=\gamma_{0}
\|\Sigma^{-1}_{n}\|\leq C ,
\end{eqnarray}
and
\begin{eqnarray}
\|\mathbf{D}^{-1}_{n}(q)-\mathbf{D}^{-1}_{n}\| \leq
\|\mathbf{D}^{-1}_{n}\|\|\mathbf{D}^{-1}_{n}(q)\|
\|\mathbf{D}_{n}(q)-\mathbf{D}_{n}\|\leq C (\sigma^{2}_{q}-\sigma^{2}_{\alpha})
\leq C\sum_{j \geq q+1}a^{2}_{j}.
\end{eqnarray}
According to (A.21)-(A.25), it remains to prove that
\begin{eqnarray}
\|\mathbf{T}_{n}(q)-\mathbf{T}_{n}\|^{2} \leq C
\sum_{j\geq q+1}|a_{j}|\sum_{j\geq q+1}j|a_{j}|.
\end{eqnarray}
By making use of Theorem 2.2 of Baxter (1962), it can be shown that
\begin{eqnarray}
\|\mathbf{T}_{n}(q)-\mathbf{T}_{n}\|_{\infty} \leq C
\sum_{j\geq q+1}|a_{j}|,
\end{eqnarray}
and
\begin{eqnarray}
\|\mathbf{T}_{n}(q)-\mathbf{T}_{n}\|_{1} \leq C
\sum_{j\geq q+1}j|a_{j}|,
\end{eqnarray}
where for an $m \times n$ matrix $\mathbf{B}=(b_{ij})_{1\leq i \leq m, 1\leq j \leq n}$,
$\|\mathbf{B}\|_{1}=\max_{1\leq j\leq n}\sum_{i=1}^{m}|b_{ij}|$
and
$\|\mathbf{B}\|_{\infty}=\max_{1\leq i\leq m}\sum_{j=1}^{n}|b_{ij}|$.
The desired conclusion (A.26) now follows from (A.27), (A.28)
and $\|\mathbf{T}_{n}(q)-\mathbf{T}_{n}\|^{2}\leq \|\mathbf{T}_{n}(q)-\mathbf{T}_{n}\|_{1}\|\mathbf{T}_{n}(q)-\mathbf{T}_{n}\|_{\infty}$.

\vspace{0.5cm}
\noindent
{\it Proof of Theorem 3.}
We will first show that for each $1\leq k \leq q_{n}$
\begin{eqnarray}
\mathrm{E}
\left\| \frac{1}{N_{0}} \sum_{t=q_{n}}^{n-1} \hat{\mathbf{e}}_{t} ( k ) \hat{e}_{t+1, k}
\right\|^{r} \leq C \left( \frac{k}{n} \right)^{r/2},
\end{eqnarray}
where $N_{0} = n - q_{n}$,
$ \hat{\mathbf{e}}_{t} ( k )=(\hat{e}_{t}, \ldots, \hat{e}_{t-k+1})^{'}$
and $\hat{e}_{t+1, k}=\hat{e}_{t+1}+\mathbf{a}^{'}(k)\hat{\mathbf{e}}_{t} ( k )$
with $\mathbf{a}^{'}(k)=(a_{1}(k), \ldots, a_{k}(k))$.
Define
$\mathbf{e}_{t} ( k ) = ( e_{t}, \ldots, e_{t-k+1} )', e_{t+1, k} = e_{t+1} +
\mathbf{a}' ( k ) \mathbf{e}_{t} ( k )$ and
$\mathbf{z}_{n}=(z_{1}, \ldots, z_{n})^{'}
= ( I - H_{d_{n}} ) \mathbf{w}(d_{n})$,
where $\mathbf{w}(d_{n})=(w_{1}(d_{n}), \ldots, w_{n}(d_{n}))^{'}=(\sum_{j=d_{n}}^{\infty}\theta_{j}x_{1j}, \ldots,
\sum_{j=d_{n}+1}^{\infty}\theta_{j}x_{nj})^{'}$.
Let $\{\mathbf{o}_{i}=(o_{1i}, \ldots, o_{ni})^{'}, 1\leq i \leq d_{n}\}$
be an orthonormal basis of the column space of $X(d_{n})$.
Then, it holds that
$H_{d_{n}}=\sum_{i=1}^{d_{n}}\mathbf{o}_{i}\mathbf{o}^{'}_{i}$.
Moreover, one has for $1\leq k \leq q_{n}$,
\begin{eqnarray}
&&\frac{1}{N_{0}} \sum_{t=q_{n}}^{n-1} \hat{\mathbf{e}}_{t} ( k ) \hat{e}_{t+1, k} \nonumber \\
&=&
\frac{1}{N_{0}} \sum_{t=q_{n}}^{n-1} \left\{ \mathbf{e}_{t} ( k ) + \mathbf{z}_{t} ( k ) -
\sum_{i=1}^{d_{n}} v_{i} \mathbf{o}_{t}^{(i)} ( k )  \right\}
\left\{ e_{t+1, k} + z_{t+1, k} -
\sum_{i=1}^{d_{n}} v_{i} o_{t+1, k}^{(i)}  \right\} \nonumber \\
&=& \displaystyle{\frac{1}{N_{0}}} \displaystyle{\sum_{t=q_{n}}^{n-1}} \mathbf{e}_{t} ( k ) e_{t+1, k} +
\displaystyle{\frac{1}{N_{0}}} \displaystyle{\sum_{t=q_{n}}^{n-1}} \mathbf{z}_{t} ( k ) e_{t+1, k} -
\displaystyle{\sum_{i=1}^{d_{n}}} v_{i} \left\{ \displaystyle{\frac{1}{N_{0}}} \displaystyle{\sum_{t=q_{n}}^{n-1}}
\mathbf{o}_{t}^{(i)} ( k ) e_{t+1, k} \right\} \nonumber \\
&+& \displaystyle{\frac{1}{N_{0}}} \displaystyle{\sum_{t=q_{n}}^{n-1}} \mathbf{e}_{t} ( k ) z_{t+1, k} +
\displaystyle{\frac{1}{N_{0}}} \displaystyle{\sum_{t=q_{n}}^{n-1}} \mathbf{z}_{t} ( k ) z_{t+1, k} -
\displaystyle{\sum_{i=1}^{d_{n}}} v_{i} \left\{ \displaystyle{\frac{1}{N_{0}}} \displaystyle{\sum_{t=q_{n}}^{n-1}}
\mathbf{o}_{t}^{(i)} ( k ) z_{t+1, k}  \right\} \nonumber \\
&-& \displaystyle{\sum_{i=1}^{d_{n}}} v_{i} \left\{ \displaystyle{\frac{1}{N_{0}}} \displaystyle{\sum_{t=q_{n}}^{n-1}}
\mathbf{e}_{t} ( k ) o_{t+1, k}^{(i)} \right\}  -
\displaystyle{\sum_{i=1}^{d_{n}}} v_{i} \left\{ \displaystyle{\frac{1}{N_{0}}} \displaystyle{\sum_{t=q_{n}}^{n-1}}
\mathbf{z}_{t} ( k ) o_{t+1, k}^{(i)} \right\} \nonumber \\
&+& \displaystyle{\sum_{i=1}^{d_{n}}} \displaystyle{\sum_{j=1}^{d_{n}}}
v_{i} v_{j} \left\{ \displaystyle{\frac{1}{N_{0}}} \displaystyle{\sum_{t=q_{n}}^{n-1}}
\mathbf{o}_{t}^{(i)} ( k ) o_{t+1, k}^{(j)} \right\} \nonumber \\ &:=& {\rm (I)} + \cdots+ {\rm (IX)},
\end{eqnarray}
where
$\mathbf{z}_{t} ( k ) = ( z_{t}, \ldots,  z_{t - k + 1} )',
z_{t+1, k} = z_{t+1} + \mathbf{a}' ( k ) \mathbf{z}_{t} ( k ),
v_{i} = \mathbf{o}_{i}' \mathbf{e},
\mathbf{o}_{t}^{(i)} ( k ) = ( o_{ti}, \ldots, o_{t - k + 1, i} )'$
and $o_{t+1, k}^{(i)} =
o_{t+1, i} + \mathbf{a}' ( k ) \mathbf{o}_{t}^{(i)} ( k )$.
By Lemmas 3 and 4 of Ing and Wei (2003),
\begin{eqnarray}
\mathrm{E} \| {\rm ( I )} \|^{r} \leq C \left( \displaystyle{\frac{k}{n}} \right)^{r/2}.
\end{eqnarray}
Theorem 2.2 of Baxter (1962) and (3.6) ensure that the spectral density of $e_{t+1, k}$ is bounded above,
and hence by (3.5), Lemma 2 of Wei (1987) and Minkowski's Inequality,
\begin{eqnarray*}
\mathrm{E} \left| \displaystyle{\frac{1}{N_{0}}} \displaystyle{\sum_{t=q_{n}}^{n-1}} z_{t-l} e_{t+1, k} \right|^{r} &\leq&
C n^{- r} \mathrm{E} \left( \displaystyle{\sum_{t=q_{n}}^{n-1}} z_{t-l}^{2} \right)^{r/2} \\
&\leq& C n^{- r/2} \mathrm{E} \left( \displaystyle{\frac{1}{n}} \displaystyle{\sum_{t=1}^{n}} w^{2}_{t}(d_{n})
\right)^{r/2} \leq C n^{- r/2} \left( \displaystyle{\sum_{j > d_{_{n}}}} | \theta_{j} | \right)^{r},
\end{eqnarray*}
for all $0 \leq l \leq k - 1$.
As a result,
\begin{eqnarray}
\mathrm{E} \| {\rm ( II )} \|^{r} \leq C \left( \displaystyle{\frac{k}{n}} \right)^{r/2}
\left( \displaystyle{\sum_{j > d_{_{n}}}} | \theta_{j} | \right)^{r}.
\end{eqnarray}
By (A.13), (3.5), the boundedness of the spectral density of $e_{t+1, k}$, and Lemma 2 of Wei (1987), one has for all $1\leq i \leq d_{n}$
and all $0 \leq l \leq k - 1$,
\begin{eqnarray}
\mathrm{E}|v_{i}|^{2r} \leq C \mathrm{E}(\sum_{t=1}^{n}o^{2}_{ti})^{r} \leq C,
\end{eqnarray}
and
\begin{eqnarray}
\mathrm{E} \left| \displaystyle{\frac{1}{N_{0}}} \displaystyle{\sum_{t=q_{n}}^{n-1}} o_{t-l, i} e_{t+1, k} \right|^{2r}
\leq C n^{- 2 r} \mathrm{E} \left( \displaystyle{\sum_{t=1}^{n}} o_{ti}^{2} \right)^{r} \leq C n^{- 2 r}.
\end{eqnarray}
Making use of (A.33), (A.34), the convexity of $x^{r}, x \geq 0$, and the Cauchy-Schwarz inequality, we obtain
\begin{eqnarray}
\mathrm{E} \| {\rm ( III )} \|^{r} &\leq& \mathrm{E} \left( \displaystyle{\sum_{i=1}^{d_{n}}} | v_{i} | \left\| \displaystyle{\frac{1}{N_{0}}}
\displaystyle{\sum_{t=q_{n}}^{n-1}} \mathbf{o}_{t}^{(i)} ( k ) e_{t+1, k} \right\| \right)^{r} \nonumber \\
&\leq& d_{n}^{r} d_{n}^{-1} \displaystyle{\sum_{i=1}^{d_{n}}} \mathrm{E} \left( | v_{i} |^{r} \left\| \displaystyle{\frac{1}{N_{0}}}
\displaystyle{\sum_{t=q_{n}}^{n-1}} \mathbf{o}_{t}^{(i)} ( k ) e_{t+1, k} \right\|^{r} \right) \nonumber \\
&\leq& C d_{n}^{r} d_{n}^{-1} \displaystyle{\sum_{i=1}^{d_{n}}} \left( \mathrm{E} \left\| \displaystyle{\frac{1}{N_{0}}}
\displaystyle{\sum_{t=q_{n}}^{n-1}} \mathbf{o}_{t}^{(i)} ( k ) e_{t+1, k} \right\|^{2r} \right)^{1/2} \nonumber \\
&\leq& C \displaystyle{\frac{d_{n}^{r} k^{r/2}}{n^{r}}}.
\end{eqnarray}
Following an argument similar to that used to prove (A.32), we have for $0 \leq l \leq k - 1$,
\begin{eqnarray*}
\mathrm{E} \left| \displaystyle{\frac{1}{N_{0}}} \displaystyle{\sum_{t=q_{n}}^{n-1}} z_{t+1, k} e_{t-l} \right|^{r}
&\leq& C n^{- r} \mathrm{E} \left( \displaystyle{\sum_{t=q_{n}}^{n-1}} z_{t+1, k}^{2} \right)^{r/2} \\
&\leq& C n^{- r} \mathrm{E} \left\{ ( 1 + \sum_{j=1}^{k}| a_{j} ( k ) | )^{2} \left( \displaystyle{\sum_{t=1}^{n}}
w_{t}^{2}(d_{n}) \right) \right\}^{r/2} \\ &\leq& C n^{- r/2} \left( \displaystyle{\sum_{j > d_{n}}}
| \theta_{j} | \right)^{r},
\end{eqnarray*}
and hence
\begin{eqnarray}
\mathrm{E} \| {\rm ( IV )} \|^{r} \leq C \left( \displaystyle{\frac{k}{n}} \right)^{r/2}  \left( \displaystyle{\sum_{j > d_{n}}}
| \theta_{j} | \right)^{r}.
\end{eqnarray}
Similarly, for all $0 \leq l \leq k - 1$,
\begin{eqnarray*}
\mathrm{E} \left| \displaystyle{\frac{1}{N_{0}}} \displaystyle{\sum_{t=q_{n}}^{n-1}} z_{t-l} z_{t+1, k} \right|^{r} \leq
C \mathrm{E} \left( \displaystyle{\frac{1}{N_{0}}} \displaystyle{\sum_{t=1}^{n}} w_{t}^{2}(d_{n}) \right)^{r} \leq
C \left( \displaystyle{\sum_{j > d_{n}}} | \theta_{j} | \right)^{2 r},
\end{eqnarray*}
yielding
\begin{eqnarray}
\mathrm{E} \| {\rm ( V )} \|^{r} \leq C \left( \displaystyle{\frac{k}{n}} \right)^{r/2}  \left( n^{1/4} \displaystyle{\sum_{j > d_{n}}}
| \theta_{j} | \right)^{2 r}.
\end{eqnarray}
By an argument analogous to (A.35), it holds that
\begin{eqnarray}
\mathrm{E} \| {\rm ( VII )} \|^{r} \leq C \displaystyle{\frac{d_{n}^{r} k^{r/2}}{n^{r}}}.
\end{eqnarray}

According to (A.33),
\begin{eqnarray*}
\mathrm{E} \| {\rm ( VI )} \|^{r} &\leq& C d_{n}^{r} d_{n}^{-1} \displaystyle{\sum_{i=1}^{d_{n}}} \left\{\mathrm{E} | v_{i} |^{r}
\mathrm{E} \left\| \displaystyle{\frac{1}{N_{0}}} \displaystyle{\sum_{t=q_{n}}^{n-1}} \mathbf{o}_{t}^{(i)} ( k ) z_{t+1, k} \right\|^{r}\right\}
\\ &\leq& C d_{n}^{r} d_{n}^{-1} \displaystyle{\sum_{i=1}^{d_{n}}}
E \left\| \displaystyle{\frac{1}{N_{0}}} \displaystyle{\sum_{t=q_{n}}^{n-1}} \mathbf{o}_{t}^{(i)} ( k ) z_{t+1, k} \right\|^{r}.
\end{eqnarray*}
Moreover, Minkowski's inequality and the Cauchy-Schwarz inequality yield
that for all $1 \leq i \leq d_{n}$ and $0 \leq l \leq k - 1$,
\begin{eqnarray*}
\mathrm{E} \left| \displaystyle{\frac{1}{N_{0}}} \displaystyle{\sum_{t=q_{n}}^{n-1}} o_{t-l, i} z_{t+1, k} \right|^{r} &\leq&
\mathrm{E} \left\{ \left( \displaystyle{\frac{1}{N_{0}}}
\displaystyle{\sum_{t=q_{n}}^{n-1}} o_{t-l, i}^{2} \right)^{r/2}
\left( \displaystyle{\frac{1}{N_{0}}} \displaystyle{\sum_{t=q_{n}}^{n-1}} Z_{t+1, k}^{2} \right)^{r/2} \right\} \\
&\leq& C n^{- r/2} \mathrm{E} \left( \displaystyle{\frac{1}{N_{0}}} \displaystyle{\sum_{t=q_{n}}^{n-1}} z_{t+1, k}^{2} \right)^{r/2}
\leq C n^{- r/2} \left( \displaystyle{\sum_{j > d_{n}}} | \theta_{j} | \right)^{r}
\end{eqnarray*}
As a result,
\begin{eqnarray}
\mathrm{E} \| {\rm ( VI )} \|^{r} \leq C \displaystyle{\frac{k^{r/2}}{n^{r/2}}} \left( d_{n}
\displaystyle{\sum_{j > d_{n}}} | \theta_{j} | \right)^{r}.
\end{eqnarray}
Similarly, it can be shown that
\begin{eqnarray}
\mathrm{E} \| {\rm ( VIII )} \|^{r} \leq C \displaystyle{\frac{k^{r/2}}{n^{r/2}}} \left( d_{n}
\displaystyle{\sum_{j > d_{n}}} | \theta_{j} | \right)^{r},
\end{eqnarray}
and
\begin{eqnarray}
\mathrm{E} \| {\rm ( IX )} \|^{r} \leq C \displaystyle{\frac{d_{n}^{2 r} k^{r/2}}{n^{r}}}.
\end{eqnarray}
Consequently, (A.29) follows from (A.30)-(A.32)
and (A.35)-(A.41).

Let
\begin{eqnarray*}
G_{n} = \displaystyle{\max_{1 \leq k \leq q_{n}}} \left\| \displaystyle{\frac{1}{N_{0}}} \displaystyle{\sum_{t=q_{n}}^{n-1}}
\hat{\mathbf{e}}_{t} ( k ) \hat{e}_{t+1, k} \right\|.
\end{eqnarray*}
Then, for any $M>0$, one obtains from (A.29) and Chebyshev's inequality that
\begin{eqnarray*}
P \left( G_{n} > M \displaystyle{\frac{q_{n}^{1/2 + 1/r}}{n^{1/2}}} \right) &=& P \left( G^{r}_{n} > M^{r}
\left( \displaystyle{\frac{q_{n}}{n}} \right)^{r/2} q_{n} \right) \\ &\leq&
\displaystyle{\frac{C}{q_{n}^{r/2 + 1} M^{r}}} \displaystyle{\sum_{k=1}^{q_{n}}} k^{r/2} \leq
\displaystyle{\frac{C}{M^{r}}}.
\end{eqnarray*}
Hence
\begin{eqnarray}
\displaystyle{\max_{1 \leq k \leq q_{n}}} \left\| \displaystyle{\frac{1}{N_{0}}} \displaystyle{\sum_{t=q_{n}}^{n-1}}
\hat{\mathbf{e}}_{t} ( k ) \hat{e}_{t+1, k} \right\| = O_{p} \left( \displaystyle{\frac{q_{n}^{1/2 + 1/r}}{n^{1/2}}} \right).
\end{eqnarray}

\vspace{0.3cm}
In the following, we shall show that
\begin{eqnarray}
\left\| \displaystyle{\frac{1}{N_{0}}} \displaystyle{\sum_{t=q_{n}}^{n-1}}
\hat{\mathbf{e}}_{t} ( q_{n} ) \hat{\mathbf{e}}_{t}' ( q_{n} )
- \Sigma_{q_{n}}  \right\| = o_{p} ( 1 ).
\end{eqnarray}
Note first that
\begin{eqnarray*}
& & \left\| \displaystyle{\frac{1}{N_{0}}} \displaystyle{\sum_{t=q_{n}}^{n-1}} \hat{\mathbf{e}}_{t} ( q_{n} )
\hat{\mathbf{e}}_{t}' ( q_{n} ) - \displaystyle{\frac{1}{N_{0}}} \displaystyle{\sum_{t=q_{n}}^{n-1}} \mathbf{e}_{t} ( q_{n} )
\mathbf{e}_{t}' ( q_{n} ) \right\| \nonumber \\
&\leq& C \left\{ \displaystyle{\frac{1}{N_{0}}} \displaystyle{\sum_{t=q_{n}}^{n-1}}
\left\| \hat{\mathbf{e}}_{t} ( q_{n} ) - \mathbf{e}_{t} ( q_{n} ) \right\|^{2} \right. \nonumber \\ &+& \left. \left( \displaystyle{\frac{1}{N_{0}}}
\displaystyle{\sum_{t=q_{n}}^{n-1}} \left\| \hat{\mathbf{e}}_{t} ( q_{n} ) - \mathbf{e}_{t} ( q_{n} ) \right\|^{2} \right)^{1/2}
\left( \displaystyle{\frac{1}{N_{0}}} \displaystyle{\sum_{t=q_{n}}^{n-1}} \left\| \mathbf{e}_{t} ( q_{n} ) \right\|^{2} \right)^{1/2}
\right\}.
\end{eqnarray*}
Straightforward calculations imply
\begin{eqnarray*}
\displaystyle{\frac{1}{N_{0}}} \displaystyle{\sum_{t=q_{n}}^{n-1}} \left\| \hat{\mathbf{e}}_{t} ( q_{n} ) - \mathbf{e}_{t} ( q_{n} )
\right\|^{2} \leq C \displaystyle{\frac{q_{n}}{N_{0}}} \| \hat{\mathbf{e}_{n}} - \mathbf{e}_{n} \|^{2} \leq  \displaystyle{\frac{Cq_{n}}{N_{0}}}
\left\{ \mathbf{e}_{n}^{'} H_{d_{n}} \mathbf{e}_{n} + \| \mathbf{w}_{n}(d_{n}) \|^{2} \right\},
\end{eqnarray*}
$\mathrm{E} ( \mathbf{e}^{'}_{n} H_{d_{n}} \mathbf{e}_{n} ) \leq C d_{n}$, $\mathrm{E} (
\| \mathbf{w}_{n}(d_{n}) \|^{2} ) \leq nC ( \sum_{j > d_{n}} | \theta_{j} | )^{2}$,
and $\mathrm{E}
( N_{0}^{-1} \sum_{t=q_{n}}^{n-1} \| \mathbf{e}_{t} ( q_{n} ) \|^{2} ) \leq C q_{n}$.
As a result,
\begin{eqnarray}
& & \left\| \displaystyle{\frac{1}{N_{0}}} \displaystyle{\sum_{t=q_{n}}^{n-1}} \hat{\mathbf{e}}_{t} ( q_{n} ) \hat{\mathbf{e}}_{t}' ( q_{n} )
- \displaystyle{\frac{1}{N_{0}}} \displaystyle{\sum_{t=q_{n}}^{n-1}} \mathbf{e}_{t} ( q_{n} ) \mathbf{e}_{t}' ( q_{n} ) \right\| \nonumber \\
&=& O_{p} \left( \displaystyle{\frac{q_{n} d_{n}}{n}} + q_{n} \left( \displaystyle{\sum_{j > d_{n}}} | \theta_{j} | \right)^{2}
+ \displaystyle{\frac{q_{n} d_{n}^{1/2}}{n^{1/2}}} + q_{n} \displaystyle{\sum_{j > d_{n}}} | \theta_{j} | \right) = o_{p} ( 1 ).
\end{eqnarray}
Moreover, Lemma 2 of Ing and Wei (2003) yields
\begin{eqnarray}
\left\| \displaystyle{\frac{1}{N_{0}}} \displaystyle{\sum_{t=q_{n}}^{n-1}} \mathbf{e}_{t} ( q_{n} ) \mathbf{e}_{t}' ( q_{n} )
- \Sigma_{q_{n}}  \right\| = O_{p} \left( \displaystyle{\frac{q_{n}}{n^{1/2}}} \right) = o_{p} ( 1 ).
\end{eqnarray}
Combining (A.44) and (A.45) leads to the desired conclusion (A.43).

\vspace{0.3cm}
By making use of (A.42) and (A.43), we next show that
\begin{eqnarray}
\left\| \hat{\mathbf{T}}_{n} ( q_{n} ) - \mathbf{T}_{n} ( q_{n} ) \right\| =& O_{p} \left(
\displaystyle{\frac{q_{n}^{1 + 1/r}}{n^{1/2}}} \right).
\end{eqnarray}

It follows from (A.43) and (A.13) that
\begin{eqnarray}
\lim_{n \to \infty}P(Q_{n}) \equiv \lim_{n \to \infty}
P(( N_{0}^{-1} \sum_{t=q_{n}}^{n-1} \hat{\mathbf{e}}_{t} ( q_{n} ) \hat{\mathbf{e}}_{t}' ( q_{n} )
 )^{-1} \,\, \mbox{exists})=1,
\end{eqnarray}
and
\begin{eqnarray}
\left\| \left( \displaystyle{\frac{1}{N_{0}}} \displaystyle{\sum_{t=q_{n}}^{n-1}} \hat{\mathbf{e}}_{t} ( q_{n} ) \hat{\mathbf{e}}_{t}' ( q_{n} )
\right)^{-1} \right\|I_{Q_{n}} = O_{p} ( 1 ).
\end{eqnarray}
In addition, an argument given in Proposition 3.1 of Ing, Chiou and Guo (2013) implies
\begin{eqnarray}
\left\| \hat{\mathbf{T}}_{n} ( q_{n} ) - \mathbf{T}_{n} ( q_{n} ) \right\|
\leq Cq^{1/2}_{n} \max_{1\leq k \leq q_{n}}\|\hat{\mathbf{a}}(k)-\mathbf{a}(k)\|,
\end{eqnarray}
where $\hat{\mathbf{a}}(k)=(\hat{a}_{1}(k), \ldots, a_{k}(k))^{'}$.
Since on $Q_{n}$,
\begin{eqnarray}
\max_{1\leq k \leq q_{n}}\|\hat{\mathbf{a}}(k)-\mathbf{a}(k)\| \leq
\left\|\left( \displaystyle{\frac{1}{N_{0}}} \displaystyle{\sum_{t=q_{n}}^{n-1}} \hat{\mathbf{e}}_{t} ( q_{n} ) \hat{\mathbf{e}}_{t}' ( q_{n} )
\right)^{-1}\right\|G_{n},
\end{eqnarray}
(A.46) is ensured by (A.42) and (A.47)-(A.50).

\vspace{0.3cm}
The proof of (3.11) is also reliant on
\begin{eqnarray}
\|\hat{\mathbf{D}}^{-1}_{n}(q_{n})-\mathbf{D}^{-1}_{n}(q_{n})\|=O_{p}\left(\frac{q^{1/r}_{n}}{n^{1/2}}+\frac{q_{n}^{1+(2/r)}}{n}\right),
\end{eqnarray}
which is in turn implied by (A.23) and
\begin{eqnarray}
\|\hat{\mathbf{D}}_{n}(q_{n})-\mathbf{D}_{n}(q_{n})\|=O_{p}\left(\frac{q^{1/r}_{n}}{n^{1/2}}+\frac{q_{n}^{1+(2/r)}}{n}\right).
\end{eqnarray}
To prove (A.52), note first that on the set $Q_{n}$,
\begin{eqnarray}
&&\max_{1 \leq k \leq q_{n}}|\hat{\sigma}^{2}_{k}-\sigma^{2}_{k}| \nonumber \\
&\leq&
\displaystyle{\max_{1 \leq k \leq q_{n}}} \left| \displaystyle{\frac{1}{N_{0}}} \displaystyle{\sum_{t=q_{n}}^{n-1}}
\hat{e}_{t+1, k}^{2} - \sigma_{k}^{2} \right| +
\left\| \left(\displaystyle{\frac{1}{N_{0}}}\displaystyle{\sum_{t=q_{n}}^{n-1}}
\hat{\mathbf{e}}_{t} ( q_{n} )\hat{\mathbf{e}}^{'}_{t} ( q_{n} )\right)^{-1}
 \right\|G^{2}_{n}.
\end{eqnarray}
Moreover, by (3.8), (3.9), Lemma 6 of Ing and Wei (2005) and an argument
similar to that used to prove (A.29), it holds that for all $1\leq k \leq q_{n}$,
\begin{eqnarray*}
\mathrm{E} \left| \displaystyle{\frac{1}{N_{0}}} \displaystyle{\sum_{t=q_{n}}^{n-1}}
\hat{e}_{t+1, k}^{2} - \sigma_{k}^{2} \right|^{r} \leq Cn^{-r/2},
\end{eqnarray*}
and hence
\begin{eqnarray}
\displaystyle{\max_{1 \leq k \leq q_{n}}} \left| \displaystyle{\frac{1}{N_{0}}} \displaystyle{\sum_{t=q_{n}}^{n-1}}
\hat{e}_{t+1, k}^{2} - \sigma_{k}^{2} \right|=O_{p}\left(\frac{q^{1/r}_{n}}{n^{1/2}}\right).
\end{eqnarray}
Similarly, we have
\begin{eqnarray}
|\hat{\gamma}_{0}-\gamma_{0}|=O_{p}(n^{-1/2}).
\end{eqnarray}
Combining (3.10), (A.42), (A.47), (A.48), (A.53)-(A.55) and
\begin{eqnarray*}
\|\hat{\mathbf{D}}_{n}(q_{n})-\mathbf{D}_{n}(q_{n})\| =
\max\{|\hat{\gamma}_{0}-\gamma_{0}|, \max_{1 \leq k \leq q_{n}}|\hat{\sigma}^{2}_{k}-\sigma^{2}_{k}|\}
\end{eqnarray*}
yields the desired conclusion (A.52).
The proof is completed by noticing
that (3.11) is an immediate consequence of (A.20), (A.46) and (A.51).

\vspace{0.5cm}
\centerline{\bf\large ACKNOWLEDGMENTS}

\vspace{0.2cm}
The research of Ching-Kang Ing was supported in part by
the Academia Sinica Investigator Award, and that of Shu-Hui Yu was
partially supported by the National Science Council of Taiwan under grant NSC 99-2118-M-390-002.
We would like to thank the editors and two
anonymous referees for their insightful and constructive comments, which greatly
improve the presentation of this paper.

\vspace{0.5cm}
\centerline{\bf\large REFERENCES}

\begin{itemize}

\item[\mbox{}] H. Akaike (1974).
{\it A new look at the statistical model identification.} \textit{IEEE Trans. Automatic Control} \textbf{19} 716-723.

\item[\mbox{}] T. Ando and K.-C. Li (2014). A model-averaging approach for high-dimensional regression. \textit{J. Amer. Statist. Assoc.} forthcoming.

\item[\mbox{}] D. W. K. Andrews (1991). Asymptotic optimality of generalized $C_{L}$, cross-validation, and generalized cross-validation in regression with heteroskedastic errors. \textit{J. Economet.} \textbf{4} 359--377.

\item[\mbox{}] G. Baxter (1962).  An Asymptotic Result for the Finite Predictor. \textit{Math. Scand.} \textbf{10} 137--144.

\item[\mbox{}]  D. R. Brillinger (1975). \textit{Time Series: Data Analysis and Theory}. Holt, Rinehart and
Winston, New York.





\item[\mbox{}] N. H. Chan, S.-F. Huang and C.-K. Ing (2013). Moment bound and mean squared prediction errors of long-memory time series. \textit{Ann. Statist.} \textbf{41} 1268--1298.


\item[\mbox{}] D. F. Findley and C. Z. Wei (1993). Moment bounds for deriving time series CLT's and model selection procedures. {\it Statist. Sinica} {\bf 3} 453--470.




\item[\mbox{}] B. E. Hansen (2007). Least squares model averaging.
\textit{Econometrica} \textbf{75} 1175--1189.


\item[\mbox{}] B. E. Hansen and J. S. Racine (2012). Jacknife model averaging. \textit{J. Economet.} \textbf{167} 38--46.



\item[\mbox{}] C.-K. Ing (2007). Accumulated prediction errors, information criteria and optimal forecasting for autoregressive time series. \textit{Ann. Statist.} {\bf 35}  1238--1277.

\item[\mbox{}] C.-K. Ing, H. T. Chiou and M. H. Guo (2013).
Estimation of inverse autocovariance matrices for long memory processes. Technical Report.

\item[\mbox{}]  C.-K. Ing and  T. L. Lai (2011). A stepwise regression method and consistent model selection for high-dimensional sparse linear models. {\it Statist. Sinica}  {\bf 21} 1473--1513.

\item[\mbox{}]  C.-K. Ing and  C.-Z. Wei (2003). On same-realization prediction in an infinite-order
autoregressive process.   {\it J.  Multivariate Anal.}  {\bf 85} 130--155.

\item[\mbox{}] C.-K. Ing and C.-Z. Wei (2005). Order selection for same-realization predictions in autoregressive processes.  \textit{Ann. Statist.} {\bf 33}  2423--2474.





\item[\mbox{}] T. L. Lai and C.-Z. Wei (1982).
Least squares estimates in stochastic regression models with applications to identification and control of dynamic systems.
\textit{Ann. Statist.} \textbf{1} 154-166.

\item[\mbox{}] S. Lee and A. Karagrigoriou (2001).  An asymptotically optimal selection of the order of a linear process. \textit{Sankhya A} \textbf{63} 93--106.




\item[\mbox{}] G. Leung and A. R. Barron (2006).
Information theory and mixing least-squares regressions.
\textit{IEEE Trans. Information Theory} \textbf{52} 3396--3410.

\item[\mbox{}] K.-C. Li (1987). Asymptotic optimality for $C_{p}$, $C_{L}$, cross-validation and generalized crossvalidation:
discrete index set. \textit{Ann. Statist.} \textbf{15} 958--975.


\item[\mbox{}] Q. Liu and R. Okui (2013). Heteroscedasticity-robust $C_{p}$ Model Averaging. \textit{Economet. J.} \textbf{16} 463-472.

\item[\mbox{}] Q. Liu, R. Okui, and A. Yoshimura (2013). Generalized Least Squares Model Averaging. Technical Report.


\item[\mbox{}] C. L. Mallows (1973). Some comments on $C_{p}$. \textit{Technometrics} \textbf{15} 661--675.

\item[\mbox{}] T. L. McMurry and D. N. Politis (2010). Banded and tapered estimates for autocovariance matrices and the linear process bootstrap. \textit{J. Time Series Anal.} \textbf{31} 471--482.



\item[\mbox{}] J. Rissanen (1986). Order estimation by accumulated prediction errors. In \textit{Essays in Time Series and Allied Processes} (J. Gani and M. P. Priestley, eds.) \textit{J. Appl. Probab.} {\bf 23A}  55--61.

\item[\mbox{}] R. Shibata (1981). An optimal selection of regression variables. \textit{Biometrika} {\bf 68}  45--54.


\item[\mbox{}] V. N. Temlyakov (2000).
Weak greedy algorithms.
\textit{Adv. Comput. Math.}
\textbf{12} 213-227.

\item[\mbox{}] A. T. K. Wan, X. Zhang, and G. Zou (2010). Least squares model averaging by Mallows criterion. \textit{J. Economet.} \textbf{156} 227--283.

\item[\mbox{}]
Z. Wang, S. Paterlini, F. Gao, and Y. Yang (2014).
Adaptive minimax regression estimation over sparse hulls. \textit{J. Machine Learning Research} \textbf{15} 1675-1711.

\item[\mbox{}]  C. Z. Wei (1987). Adaptive prediction by least squares predictors in stochastic regression
models with applications to time series. \textit{Ann. Statist.} \textbf{15} 1667--1682.

\item[\mbox{}]
X. Wei and Y. Yang (2012).
Robust forecast combinations. \textit{J. Economet.} \textbf{166} 224-236.

\item[\mbox{}]P. Whittle (1960). Bounds for the moments of linear and quadratic forms in independent variables.
\textit{Theory Probab. Appl.} \textbf{5} 302--305.

\item[\mbox{}]W. B. Wu and M. Pourahmadi (2003). Nonparametric estimation of large covariance matrices
of longitudinal data. \textit{Biometrika} \textbf{90} 831--844.

\item[\mbox{}] W. B. Wu and M. Pourahmadi (2009). Banding sample covariance matrices of stationary processes. \textit{Statist. Sinica} \textbf{19} 1755--1768.

\item[\mbox{}] Y. Yang (2001). Adaptive regression by mixing. \textit{J. Amer. Statist. Assoc.} \textbf{96} 574--586.

\item[\mbox{}] Y. Yang (2007).
Prediction/Estimation with Simple Linear Model: Is It Really that Simple? \textit{Economet. Theory} \textbf{23} 1-36.

\item[\mbox{}] S. H. Yu, C.-C. Lin and H.-W. Cheng (2012). A note on mean squared prediction error under the unit root model with deterministic trend. \textit{J. Time Series Anal.} \textbf{33} 276--286.

\item[\mbox{}] Z. Yuan and Y. Yang (2005). Combining linear regression model: when and how? \textit{\it J. Amer. Statist. Assoc.} \textbf{100} 1202--1214.

\item[\mbox{}] X. Zhang, A. T. K. Wan and G. Zou (2013). Model averaging by jackknife criterion in models with dependent data. \textit{J. Economet.} \textbf{174} 82--94.

\item[\mbox{}] A. Zygmund (1959). Trigonometric Series, 2nd ed. Cambridge Univ. Press.

\end{itemize}

\end{document}